\newcommand{\R}{{\mathds R}}
\newcommand{\Z}{\mathds Z}

\newcommand{\N}{{\mathds N}}
\newcommand{\C}{\mathds C}
\newcommand{\Gr}{\mathrm{Gr}}

\newcommand{\Deltao}{\Delta\!^{\scriptscriptstyle o}}

\newcommand{\lip}{\left[\!\left[}
\newcommand{\rip}{\right]\!\right]}
\newcommand{\biglip}{\left[\!\!\left[}
\newcommand{\bigrip}{\right]\!\!\right]}

\newcommand{\Dim}{\mathrm{dim}}

\newcommand{\Ker}{\mathrm{Ker}}
\newcommand{\complex}[1]{\,{\vphantom{#1}}^\textrm{c}\!\!#1}
\newcommand{\ccomplex}[1]{\,{\vphantom{#1}}^{\bar{\textrm{c}}}\!#1}

\renewcommand{\contentsline}[3]{\csname new#1\endcsname{#2}{#3}}
\newcommand{\newchapter}[2]{\bigskip\hbox to \hsize{\vbox{\advance\hsize by -.5cm\baselineskip=12pt\parfillskip=0pt\leftskip=2cm\noindent\hskip -2cm #1\leaders\hbox{.}\hfil\hfil\par}$\,$#2\hfil}}
\newcommand{\newsection}[2]{\medskip\hbox to \hsize{\vbox{\advance\hsize by -.5cm\baselineskip=12pt\parfillskip=0pt\leftskip=2.5cm\noindent\hskip -2cm #1\leaders\hbox{.}\hfil\hfil\par}$\,$#2\hfil}}
\newcommand{\newsubsection}[2]{\medskip\hbox to \hsize{\vbox{\advance\hsize by -.5cm\baselineskip=12pt\parfillskip=0pt\leftskip=3.5cm\noindent\hskip -2cm #1\leaders\hbox{.}\hfil\hfil\par}$\,$#2\hfil}}

\documentclass[oneside,draft,a4paper,11pt]{amsart}

\usepackage{times} 
\usepackage{dsfont}

\numberwithin{equation}{section}

\title[Maslov index in locally symmetric semi-Riemannian manifolds]{%
Conjugate points and Maslov index in locally symmetric
semi-Riemannian manifolds.}
\author[M.\ A.\ Javaloyes]{Miguel Angel Javaloyes}
\author[P.\ Piccione]{Paolo Piccione}
\address{Departamento de Matem\'atica,\hfill\break\indent
Instituto de Matem\'atica e Estat\'\i stica\hfill\break\indent
Universidade de S\~ao Paulo, \hfill\break\indent Rua do Mat\~ao
1010, CEP 05508-900, S\~ao Paulo, SP\hfill\break\indent Brazil}
\email{majava@ime.usp.br, piccione@ime.usp.br}
\dedicatory{In memory of Enzo Baldoni, a man of peace.}


\date{May 6th, 2005}

\begin{document}


\theoremstyle{plain}\newtheorem{teo}{Theorem}[section]
\theoremstyle{plain}\newtheorem{prop}[teo]{Proposition}
\theoremstyle{plain}\newtheorem{lem}[teo]{Lemma}
\theoremstyle{plain}\newtheorem{cor}[teo]{Corollary}
\theoremstyle{definition}\newtheorem{defin}[teo]{Definition}
\theoremstyle{remark}\newtheorem{rem}[teo]{Remark}
\theoremstyle{plain} \newtheorem{assum}[teo]{Assumption}
\theoremstyle{definition}\newtheorem{example}[teo]{Example}


\begin{abstract}
We study the singularities of the exponential map in
semi Riemann\-ian locally symmetric manifolds.
Conjugate points along geodesics  depend only on
real negative eigenvalues of the curvature tensor,
and their contribution to the Maslov index of the geodesic
is computed explicitly.
We prove that degeneracy of conjugate points, which is a
phenomenon that can only occur in semi-Riemannian geometry,
is caused in the locally symmetric case by the lack of diagonalizability
of the curvature tensor. The case of Lie groups endowed with a bi-invariant
metric is studied in some detail, and conditions are given
for the lack of local injectivity of the exponential map
around its singularities.
\end{abstract}

\maketitle

\tableofcontents

\begin{section}{Introduction}\label{sec:intro}
The geodesic flow in semi-Riemannian manifolds, i.e., manifolds endowed
with a metric tensor which is not positive definite, has features which
are quite different from the Riemannian, i.e., positive definite, case.
Although the local theory of semi-Riemannian geodesics is totally equivalent
to the Riemannian one, when it gets to global properties the situation
changes dramatically. Most notably, compact manifolds may fail to
be geodesically connected, and the classical Morse theory for geodesics
does not apply to the non positive definite case. In this paper we
will be concerned with another phenomenon typical of the semi-Riemannian
world, which is the existence of degenerate singularities for the
exponential map. Unlike the Riemannian case, degenerate conjugate points
may accumulate along a geodesic, and they do not necessarily
determine bifurcation. The theoretical occurrence and the relevance
of such phenomena has been studied recently in a series of papers;
however, no explicit calculation has been carried out so far due
to the difficulties in the integration of the geodesic equation.
If one wants to study the global geometry of the conjugate
locus in a semi-Riemannian manifold, he will find somewhat
discouraging the result proven in \cite{fechado}, concerning
the distribution of conjugate points along a geodesic. Such set can
be arbitrarily complicated:  any bounded closed subset of the real
line, like Cantor sets or other pathological examples, appears as
the set of conjugate instants along spacelike geodesics in conformally
flat Lorentzian $3$-dimensional manifolds.
It is therefore hopeless to be able to develop significant results
concerning the geometry of the conjugate locus in
the general case of smooth metrics.
On the other hand,
if one restricts his attention to the case of real-analytic metrics
then accumulation does not occur, and higher order methods for
analyzing the isolated singularities of the exponential map are
available (see \cite{GiaPicPor}). As in the Riemannian case (see \cite{Mil1, Mil2}),
in order to make explicit computation, an important family of examples
of analytic semi-Riemannian manifolds to start with is given by
the class of Lie groups endowed with an invariant metric.
As a first step in this direction, in this paper we will consider
the case of (non compact) Lie groups endowed with a bi-invariant
semi-Riemannian metric or, more generally, the case of semi-Riemannian
locally symmetric spaces. Recall that if $G$ is a semi-simple Lie group,
then the Killing form of its Lie algebra $\mathfrak g$  defines a bi-invariant
semi-Riemannian metric on $G$; more generally, given  a nondegenerate symmetric
bilinear form $B$ on $\mathfrak g$ such that $\mathrm{ad}_X$ is
$B$-skew symmetric for all $X\in\mathfrak g$, then $B$ can be
extended to a bi-invariant semi-Riemannian metric on $G$.
For instance, if $G$ is semi-simple and non compact, then its Killing
form is not definite, and we obtain a non trivial class of examples
where the occurrence of several types of nondegeneracies can
be detected by explicit computations. The class of Lie groups
admitting a bi-invariant semi-Riemannian metric is quite large,
and it has been described in \cite{medina}.
In the present paper we develop an algebraic theory that allows
to determine all the singularities of the exponential
map of a locally symmetric semi-Riemannian manifold, to characterize
which of these singularities are degenerate, and we give a general
formula for computing an important integer valued invariant
for geodesics called the \emph{Maslov index}. This integer
number is given by an algebraic count of the conjugate instants
along a geodesic; the notion of Maslov index appears naturally in the
infinite dimensional Morse theory for the strongly indefinite
functionals, where it plays the role of a generalized Morse index
(see \cite{AbbBenForMas}).
\smallskip

The Riemannian curvature tensor of a locally symmetric semi-Rie\-mann\-ian
manifold $(M,g)$ is parallel, so that the Jacobi equation along a
geodesic $\gamma$ is represented, via a parallel trivialization of
the tangent bundle $TM$ along $\gamma$, by a second order linear
equation with constant coefficients. The singularities of the
exponential map of $(M,g)$ are zeroes of solutions
of such equations, and they exist when the curvature tensor
has real negative eigenvalues  (Lemma~\ref{thm:existenceconjpts}).
Degeneracies of such singularities correspond to degeneracies of
the restriction of the metric tensor $g$ to the generalized eigenspaces
of the curvature tensor relative to the real negative eigenvalues
(Proposition~\ref{thm:restrgeigen}
and Corollary~\ref{thm:cordegconjinst}).
When $G$ is a Lie group and $h$ is a bi-invariant semi-Riemannian
metric on $G$, in which case the geodesics through the identity
are the one-parameter subgroups of $G$, the conjugate points are
determined by the purely imaginary eigenvalues of the adjoint
map (Proposition~\ref{thm:gruppi}). As in the Riemannian case (see \cite{Mil1}),
the multiplicity of each conjugate point in a bi-invariant semi-Riemannian
metric is even. In the special case of a bi-invariant Lorentzian
metric on a Lie group whose dimension is less than $6$, then
the Maslov index of a geodesic equals the number of conjugate points
(counted with multiplicity) along the geodesic (Proposition~\ref{thm:gruppi2}).

The preliminary algebraic results needed to carry out our
computations are collected in Section~\ref{sec:algprel}. An effort
has been made to make the paper self-contained, and, to this aim,
in Section~\ref{sec:algprel} we have reproduced the proof of some
well known facts (see \cite{GohLanRod}) about the Jordan form of
endomorphisms that are symmetric with respect to non positive
definite inner products. New algebraic invariants called {\em
Jordan signatures\/} are introduced in
Subsection~\ref{sub:Jordansignature}; these are nonnegative
integers associated to each (real) eigenvalue of a $g$-symmetric
endomorphism, and they are used in the computation of the
contribution to the Maslov index given by the final endpoint.
Conjugate points for arbitrary differential systems are defined
and discussed in Section~\ref{sec:eigenvalues}, where we prove
that, for an arbitrary system with constant coefficients, the
conjugate instants are determined solely by the real negative
eigenvalues of the coefficient matrix.

The Maslov index is computed in Section~\ref{sec:computation} (Corollary~\ref{thm:maslovindex})
using a formula
proven in Lemma~\ref{thm:formulacompMaslovind} that relates this number with the variation
of the extended coindex of a smooth path of symmetric bilinear forms
defined on the space of Jacobi fields. Using similar formulas, another
symplectic invariant called the \emph{Conley--Zehnder index}
is computed explicitly for systems arising from the Jacobi equation
of a locally symmetric semi-Riemannian manifold.
Finally, in Section~\ref{sec:final}
we  make some explicit computations in semi-Riemannian Lie groups, and
we show how one can extend the results to more general classes
of symplectic systems.

\end{section}

\begin{section}{Algebraic preliminaries}
\label{sec:algprel} We will be concerned  with
second order linear systems whose matrix of coefficients $A$ is
symmetric relatively to a nondegenerate symmetric bilinear form
$g$ on $\R^n$, which is not necessarily positive definite. Then,
$A$ may not be diagonalizable, and in fact we will show that, when
$A$ is the curvature tensor of a semi-Riemannian metric, the
occurrence of such circumstance determines the existence of
degenerate singularities of the exponential map.

In order to study these singularities and to carry out the necessary
computations, it seems natural to use a Jordan basis for the curvature
tensor of the semi-Riemannian metric.  Following the theory in \cite{GohLanRod},
one proves that in such basis, the matrix
representation of the metric has a simple expression (see Proposition~\ref{thm:restrgeigen}),
which allows a direct computation of the Maslov index without employing
perturbation arguments. This
result will then be used to
study the restriction of $g$ to the generalized eigenspaces of $A$
and to define the notion of {\em Jordan signatures}.
\subsection{Jordan form of $g$-symmetric
endomorphisms}\label{sub:Jordanform}
Let us introduce our terminology and fix our notations by recalling
a few elementary facts concerning the Jordan
canonical form for matrices representing linear endomorphisms
of $\R^n$. Let $A:\R^n\to\R^n$ a linear
endomorphism; when needed, we will consider the $\C$-linear
extension of $\complex A$ to an endomorphism of $\C^n$, defined by $\complex A(x+iy)=Ax+iAy$.
Given a complex number $z$, we will denote by $\Im(z)$ its imaginary part.

By $\mathfrak s(A)$ we will mean the spectrum of $\complex A$; for $\lambda\in\mathfrak s(A)$,
let $\mathcal H_\lambda(A)$ denote the \emph{complex generalized eigenspace} of $A$:
\[\mathcal H_\lambda(A)=\Ker(\complex A-\lambda)^n.\]
If $\lambda\in\mathfrak s(A)$ then obviously $\overline\lambda\in\mathfrak s(A)$;
we set:
\[\mathcal F_\lambda(A)=\begin{cases}\mathcal H_\lambda(A),& \text{if}\ \lambda\in\mathfrak s(A)\cap\R;\\ \\
\mathcal H_\lambda(A)\oplus\mathcal H_{\overline\lambda}(A),&\text{if}\ \lambda\in\mathfrak s(A)\setminus\R,\end{cases}\]
so that:
\begin{equation}\label{eq:tuttoCn}
\C^n=\bigoplus_{\substack{\lambda\in\mathfrak s(A)\\ \Im(\lambda)\ge0}}\mathcal F_\lambda(A).
\end{equation}
Finally, let $\mathcal F_\lambda^o(A)$ denote the \emph{real generalized eigenspace} of $A$:
\[\mathcal F_\lambda^o(A)=\mathcal F_\lambda(A)\cap\R^n;\]
$\mathcal F_\lambda(A)$ is the complexification of $\mathcal F_\lambda^o(A)$, i.e.,
$\mathcal F_\lambda(A)=\mathcal F_\lambda^o(A)+i\mathcal F_\lambda^o(A)$, and thus
\[\R^n=\bigoplus_{\substack{\lambda\in\mathfrak s(A)\\ \Im(\lambda)\ge0}}\mathcal F^o_\lambda(A).\]
Clearly, if $\lambda\in\R$, then  $\mathrm{dim}_\C\big(\Ker(\complex A
-\lambda)\big)=\mathrm{dim}_{\R}\big(\Ker(A-\lambda)\big)$ and $\mathcal F^o_\lambda(A)=\Ker(A-\lambda)^n$;
we will call the dimension of $\Ker(A-\lambda)$  the {\em geometric multiplicity\/}
of the eigenvalue $\lambda$, while the dimension of  $\Ker(A-\lambda)^n$
will be called the {\em algebraic multiplicity} of $\lambda$.

\noindent\
The spaces $\mathcal H_\lambda(A)$ (and  $\mathcal F_\lambda(A)$) are $\complex A$-invariant,
and the restriction $\complex A\vert_{\mathcal H_\lambda(A)}$ of $\complex A$ to $\mathcal H_\lambda(A)$ is represented
in a suitable basis by a matrix which is the direct sum of \emph{$\lambda$-Jordan blocks}, i.e., matrices
of the form:
\begin{equation}\label{eq:Jordanblock}
\begin{pmatrix}\lambda&1&0&0&\ldots&0\\ 0&\lambda&1&0&\ldots&0\\
&&\vdots\\ 0&0&\ldots&\lambda&1&0\\ 0&0&\ldots&0&\lambda&1\cr 0&0&\ldots&0&0&\lambda\end{pmatrix}.
\end{equation}
By direct sum of the $k_1\times k_1$ matrix $\alpha$ and the $k_2\times k_2$ matrix
$\beta$, we mean the $(k_1+k_2)\times(k_1+k_2)$ matrix given by:
\[\alpha\oplus\beta=\begin{pmatrix}\alpha&0\\ 0&\beta\end{pmatrix}.\]
We will denote by $J_k(\lambda)$ a Jordan block of the form \eqref{eq:Jordanblock}
having size $k\times k$ when $k>1$; $J_1(\lambda)$ is defined to be the
$1\times1$ matrix $(\lambda)$.

The decomposition of $\complex A\vert_{\mathcal H_\lambda(A)}$
into direct sum of $\lambda$-Jordan blocks is not unique, but the
number of blocks (and their dimension) appearing in this
decomposition  is fixed, and it is equal to the complex dimension
of $\Ker(\complex A-\lambda)$. We will now determine the Jordan
decomposition of endomorphisms obtained from $A$ by analytic
functional calculus.
\smallskip

In what follows, we will denote by $N_r$ the $r\times r$ nilpotent matrix:
\[
N_r=\begin{pmatrix}0&1&0&0&\ldots&0\\ 0&0&1&0&\ldots&0\\
&&\vdots\\ 0&0&\ldots&0&1&0\\ 0&0&\ldots&0&0&1\cr 0&0&\ldots&0&0&0\end{pmatrix}.
\]

\begin{lem}\label{lem:seriediA}
Let $\lambda\in\C$, $B=\lambda\cdot \mathrm I_r+N_r$ and let $h:U\to\C$ be an analytic function
defined on an open  $U\subset\C$ containing $0$
whose Taylor series $h(x)=\sum_{i=0}^\infty a_i x^i$ has radius of convergence $r>\vert\lambda\vert$. Then,
$h(B)=\sum_{i=0}^\infty a_iB^i$ converges, and
\begin{equation}\label{eq:seriediA}
h(B)=\sum_{i=0}^{r-1}\frac{1}{i!}h^{(i)}(\lambda)N_r^i,
\end{equation}
where $h^{(i)}$ is the $i$-th derivative of $h$. If $h'(\lambda)\ne0$, then the canonical Jordan
form of $h(B)$ is given by:
\begin{equation*}\label{eq:Jordanblock2}
\begin{pmatrix}h(\lambda)&1&0&0&\ldots&0\\ 0&h(\lambda)&1&0&\ldots&0\\
&&\vdots\\ 0&0&\ldots&h(\lambda)&1&0\\ 0&0&\ldots&0&h(\lambda)&1\cr 0&0&\ldots&0&0&h(\lambda)\end{pmatrix}.
\end{equation*}
\end{lem}
\begin{proof}
By linearity, it suffices to prove \eqref{eq:seriediA} for the function $h(x)=x^p$, with $p\in\N$.
The proof of the desired equality in this case follows trivially from the binomial formula.
The second statement follows now easily, observing that $(h(B)-h(\lambda)\mathrm I_r)^{r-1}$
is the  matrix $h'(\lambda)^{r-1}N_r^{r-1}$.
\end{proof}
\begin{cor}\label{cor:eigenFu}
Let $A$ be an endomorphism of $\C^n$ and let $h:U\to\C$ be an analytic
function defined on an open  $U\subset\C$ containing $0$.
Assume that the Taylor series of $h$ centered at $0$ has radius of
convergence $r>\vert\lambda\vert$ for all $\lambda\in\mathfrak s(A)$;
then, $\mathfrak s(h(A))=h(\mathfrak s(A))$.
\qed
\end{cor}

Let us now consider a nondegenerate symmetric bilinear form $g(\cdot,\cdot)$ on $\R^n$;
it will be convenient to identify $g$ with the corresponding linear map\footnote{%
In this paper, the superscript ${}^*$ attached to the symbols of spaces or maps
will denote duality. When attached to matrices, it will denote the (conjugate) transpose.} $\R^n\ni v\mapsto g(\cdot,v)\in{\R^n}^*$.
Nondegeneracy means that $g$ is an isomorphism, and symmetry means that $g^*=g$.
Let $\ccomplex g$ denote the unique \emph{sesquilinear} extension of $g(\cdot,\cdot)$
to $\C^n\times\C^n$; in this case, $\ccomplex g$ will be identified with the
\emph{conjugate linear} map $\ccomplex g:\C^n\to{\C^n}^*$ obtained as the unique
conjugate linear extension of $g:\R^n\to{\R^n}^*$. Nondegeneracy of $g$ is equivalent to
the nondegeneracy of $\ccomplex g$, and the symmetry of $g$ is equivalent to
$\ccomplex g$ being conjugate symmetric, i.e., $\ccomplex g(v,w)=\overline{\ccomplex g(w,v)}$ for all
$v,w\in\C^n$.

The \emph{index $\mathrm n_-(B)$} and the \emph{coindex $\mathrm n_+(B)$}\/ of a symmetric bilinear form $B$
defined on a (finite dimensional) real vector space $V$
are defined respectively
to be the number of $-1$'s and the number of $1$'s in the canonical
matrix representation of $B$ given by Sylvester's Inertia Theorem;
by {\em signature} of $B$, denoted by $\sigma(B)$, we will mean
the difference $\mathrm n_+(B)-\mathrm n_-(B)$. The \emph{nullity} of $B$ is the
dimension of the kernel of $B$, defined by $\Ker(B)=\big\{v\in V:B(v,w)=0\ \text{for all}\ w\in V\big\}$.

A subspace $W\subset V$ is said to be \emph{$B$-positive} (resp., \emph{$B$-negative})
if $B\vert_{W}$ is\footnote{%
With a slight abuse of notations,
given a symmetric bilinear form $B$ on a vector space $V$ and given a subspace $W$ of $V$, we will
denote by $B\vert_W$ the restriction of $B$ to $W\times W$.}
positive definite (resp., negative definite); a subspace
$W\subset V$ will be called \emph{$B$-isotropic} if $B\vert_{W}$ vanishes
identically. The index (resp., the coindex) of $B$ is equal to the dimension of a maximal
$B$-negative (resp., $B$-positive) subspace of $V$.

\begin{rem}\label{thm:remsignature}
If $B$ is nondegenerate and $W\subset V$ is a {$B$-isotropic}
subspace of $V$,  then $\mathrm n_\pm(B)\ge\Dim(W)$ and
$\vert\sigma(B)\vert\le\mathrm{dim}(V)-2\,\mathrm{dim}(W)$.
Namely, if $W_-$ (resp., $W_+$) is a maximal $B$-negative (resp., $B$-positive)
subspace of $V$, then $W_\pm\cap W=\{0\}$, hence $\Dim(W_\pm)\le\Dim(V)-\Dim(W)$.
Moreover, since $B$ is nondegenerate, $\Dim(W_+)+\Dim(W_-)=\Dim(V)$, from which
the three inequalities asserted follow easily.
\end{rem}

\noindent\ \
We will assume that $A$ is $g$-symmetric, meaning that $g(Av,w)=g(v,Aw)$ for all $v,w\in\R^n$;
in terms of linear maps, this is equivalent to requiring that the following equality holds:
$gA=A^*g$. The $g$-symmetry of $A$ is equivalent to the $\ccomplex g$-simmetry
of $\complex A$.
\begin{lem}\label{thm:sommagort}
If $\lambda,\mu\in\mathfrak s(A)$ are such that $\lambda\ne\overline\mu$, then
the generalized eigenspaces $\mathcal H_\lambda(A)$ and $\mathcal H_\mu(A)$ are
$\ccomplex g$-orthogonal. If $\lambda\in\mathfrak s(A)$, then the restriction
of the bilinear form $\ccomplex g$ to $\mathcal F_\lambda(A)$ is nondegenerate, and
so is the restriction of $g$ to $\mathcal F_\lambda^o(A)$.
In particular, if $\lambda\in\mathfrak s(A)\cap\R$, then the restriction
of $g$ to $\Ker(A-\lambda)^n$ is nondegenerate.
\end{lem}
\begin{proof}
We show by induction on $k=k_1+k_2$ that
$\Ker(\complex A-\lambda)^{k_1}$ and $\Ker(\complex A-\mu)^{k_2}$ are
$\ccomplex g$-orthogonal spaces. When $k_1=k_2=1$ it is just a direct computation,
namely, for $v\in\Ker(\complex A-\lambda)$ and $w\in\Ker(\complex A-\mu)$ one has:
\begin{equation}\label{eq:vedidopo}
\lambda \ccomplex g(v,w)\!=\!\ccomplex g(\lambda v,w)\!=\!\ccomplex g(\complex A v,w)\!=\!
\ccomplex g(v,\complex Aw)\!=\!\ccomplex g(v,\mu w)\!=\!\overline\mu \ccomplex g(v,w)
\end{equation}
which implies  $\ccomplex g(v,w)=0$.

Assume now that $\Ker(\complex A-\lambda)^{k_1}$ and $\Ker(\complex A-\mu)^{k_2}$ are
$\ccomplex g$-orthogonal spaces for all pairs $k_1$ and $k_2$ such that $k_1+k_2<k$; let
$s_1,s_2\ge1$ be such that $s_1+s_2=k$, and let $v\in\Ker(\complex A-\lambda)^{s_1}$
and $w\in\Ker(\complex A-\mu)^{s_2}$. Since $(\complex A-\lambda)v\in\Ker(\complex A-\lambda)^{s_1-1}$
and $(\complex A-\mu)w\in\Ker(\complex A-\mu)^{s_2-1}$, by the induction hypothesis, we have:
\[\ccomplex g\big((\complex A-\lambda)v,w\big)=\ccomplex g\big(v,(\complex A-\mu)w\big)=0,\]
and from these two equalities it follows easily  $\ccomplex g(v,w)=0$,  as
in \eqref{eq:vedidopo}.

The orthogonality of the generalized eigenspaces shows that  \eqref{eq:tuttoCn}
is in fact a $\ccomplex g$-orthogonal direct decomposition of $\C^n$,
from which it follows that the restriction of $\ccomplex g$ to each $\mathcal F_\lambda(A)$
is nondegenerate, since $\ccomplex g$ is
nondegenerate on $\C^n$. Finally, the nondegeneracy  of the
restriction of $\ccomplex g$ on $\mathcal F_\lambda(A)$ is equivalent
to the nondegeneracy of the restriction of $g$ to $\mathcal F_\lambda^o(A)$;
in particular, if $\lambda\in\R$, then $g$ is nondegenerate on $\Ker(A-\lambda)^n$.
\end{proof}

In order to study the restriction of $g$ to the generalized eigenspaces
of $A$, we will now determine the form of the matrix representing $g$ in
a suitable Jordan basis for $A$. Lemma~\ref{thm:sommagort} tells us that it is not restrictive
to consider the case that $A$ has only two complex conjugate eigenvalues
or one real eigenvalue: once the matrix representation $g_\lambda$ of
$g\vert_{\mathcal F_\lambda^o(A)}$ has been determined for each $\lambda\in\mathfrak s(A)$
with $\Im(\lambda)\ge0$, then the matrix representation of $g$ will be given
by the direct sum of all such $g_\lambda$'s.
As a matter of facts, we will only be interested
in the case of one real eigenvalue (see Lemma~\ref{thm:existenceconjpts} below).
Using the terminology of \cite{GohLanRod}, we will call a {\em sip matrix}
an $n\times n$ matrix $\mathrm{Sip}_n$ of the form:
\begin{equation}\label{eq:sipmatrix}
\mathrm{Sip}_n=\begin{pmatrix}0&0&\ldots&0&0&1\\ 0&0&\ldots&0&1&0\\ &&&\vdots\\ 0&0&1&\ldots&0&0\\0&1&0&\ldots&0&0\\
1&0&0&\ldots&0&0\end{pmatrix}.
\end{equation}
Adapting the proof of \cite[Theorem~3.3]{GohLanRod}, we get the following:

\begin{prop}\label{thm:restrgeigen}
Let $\lambda$ be a real eigenvalue of $A$, with $r=\Dim(\Ker(A-\lambda))$.
Then, the real generalized eigenspace
can be written as a $g$-orthogonal direct sum:
\[\Ker(A-\lambda)^n=\bigoplus_{i=1}^r\mathcal V_{\lambda,i},\]
for which the following properties hold:
\begin{itemize}
\item[(a)] $g\vert_{\mathcal V_{\lambda,i}}$ is nondegenerate for all $i=1,\ldots,r$;
\item[(b)] each $\mathcal V_{\lambda,i}$ is $A$-invariant;
\item[(c)] for all $i$, there exists a basis $v_1^i,\ldots,v_{n_i}^i$ of
$\mathcal V_{\lambda,i}$ and a number $\epsilon_i\in\{-1,1\}$
such that in this basis the matrix representation of $A\vert_{\mathcal V_{\lambda,i}}$
is as in \eqref{eq:Jordanblock}, and
the matrix representation of $g\vert_{\mathcal V_{\lambda,i}}$ is given by $\epsilon_i\cdot\mathrm{Sip}_{n_i}$.
\end{itemize}
\end{prop}
\begin{proof}
It will suffice to show the existence of a number $\epsilon=\{-1,1\}$,
of a subspace $\mathcal V\subset\Ker(A-\lambda)^n$
and of a basis $w_1,\ldots,w_s$ of $\mathcal V$ with the properties:
\begin{itemize}
\item $Aw_1=\lambda w_1$ and $Aw_j=w_{j-1}+\lambda w_j$ for $j=2,\ldots,s$;
\item $g(w_j,w_k)=\epsilon\delta_{j+k,s+1}$ for all $j,k=1,\ldots,s$.
\end{itemize}
The two properties above imply that $\mathcal V$ is $A$-invariant
and that the restriction $g\vert_{\mathcal V}$ is nondegenerate. The matrix
representation of $A\vert_{\mathcal V}$  in the
basis $w_1,\!\ldots,\!w_s$ is as in \eqref{eq:Jordanblock}
and the matrix representation of $g\vert_{\mathcal V}$
is $\epsilon\cdot\mathrm{Sip}_s$; the conclusion
will follow easily from an induction argument by considering the $g$-orthogonal
complement of $\mathcal V$ in $\mathcal F_\lambda^o(A)$.

To infer the existence of such a subspace $\mathcal V$ with the desired basis,
let us argue as follows.
There exists $s\ge1$ with the property that $(A-\lambda)^s\vert_{\Ker(A-\lambda)^n}=0$
but $(A-\lambda)^{s-1}\vert_{\Ker(A-\lambda)^n}\ne0$; since $B=g\big((A-\lambda)^{s-1}\cdot,\cdot\big)$ is
a non zero symmetric bilinear form on $\Ker(A-\lambda)^n$, there must exists
a vector $a_1$ such that $B(a_1,a_1)\ne0$. We can normalize $a_1$ in such a way that
$g\big((A-\lambda)^{s-1}a_1,a_1\big)=\epsilon$, for some $\epsilon\in\{-1,1\}$; the case
$s=1$ is concluded by setting $w_1=a_1$, and we will now assume $s>1$.

For $j=1,\ldots,s$, let us define $a_j=(A-\lambda)^{j-1}a_1$  and let $\mathcal V$ be the
space spanned by the $a_j$'s; it is very easy to check that the $a_j$'s are linearly
independent, and thus $\Dim(\mathcal V)=s$. For $j+k=s+1$, we have:
\begin{equation}\label{eq:s+1}\begin{split}g(a_j,a_k)\;&=g\big((A-\lambda)^{j-1}a_1,
(A-\lambda)^{k-1}a_1\big)\\&=g\big((A-\lambda)^{j+k-2}a_1,a_1\big)=
g\big((A-\lambda)^{s-1}a_1,a_1\big)=\epsilon,\end{split}
\end{equation}
while if $j+k>s+1$ we have:
\begin{equation}\label{eq:s+2}g(a_j,a_k)=g\big((A-\lambda)^{j+k-2}a_1,a_1\big)=0.
\end{equation}
Now, set $b_1=a_1+\alpha_2 a_2+\ldots+\alpha_s a_s$ and $b_j=(A-\lambda)^{j-1}b_1$
for $j=1,\ldots,s$. Here the real coefficients $(\alpha_i)_{i=2}^s$ are to be determined
in such a way that $g(b_1,b_j)=0$ for all $j=1,\ldots,s-1$, which would imply easily $g(b_j,b_k)=\epsilon\delta_{j+k,s+1}$
for all $j$ and $k$. Such a choice of
the $\alpha_i$'s is indeed possible (and unique), namely, the equality $g(b_1,b_j)=0$  is given, in
view of \eqref{eq:s+1} and \eqref{eq:s+2}, by:
\[\begin{split}0\;&=g\big(a_1+\sum\nolimits_{k=2}^s\alpha_k a_k,a_j+\sum\nolimits_{k=2}^{s-j+1}\alpha_k a_{j+k-1}\big)\\
&=g(a_1,a_j)+2\epsilon\alpha_{s-j+1}+\text{terms in $\alpha_2,\ldots,\alpha_{s-j}$},\end{split}\]
so that the $\alpha_i$'s can be determined recursively by taking $j=s-1, s-2,\ldots,1$ in the above equality.
It is easy to check that the $b_j$'s form a basis of $\mathcal V$.
Finally, set $w_j=b_{s-j+1}$ for all $j=1,\ldots,s$;
an immediate computation shows that the $w_j$'s have the required properties.
\end{proof}
We draw a first immediate conclusion from the above result:
\begin{cor}
\label{thm:corsignature}
If $\lambda$ is a real eigenvalue of $A$, then the absolute value of the
signature of the restriction
of $g$ to $\Ker\big((A-\lambda)^n\big)$ is less than or equal to the dimension
of $\Ker(A-\lambda)$. The restriction of $g$ to the eigenspace $\Ker(A-\lambda)$ is
nondegenerate if and only if the algebraic multiplicity and the geometric multiplicity
of $\lambda$ coincide.
\end{cor}
\begin{proof}
Since the signature of $g$ is additive by $g$-orthogonal sums,
using the result of Proposition~\ref{thm:restrgeigen} it suffices
to show that $\left\vert\sigma\big(g\vert_{\mathcal V_{\lambda,i}}\big)\right\vert\le 1$
for all $i=1,\ldots,r=\Dim(\Ker(A-\lambda))$. Since $g\vert_{\mathcal V_{\lambda,i}}$
is represented by the matrix $\epsilon_i\cdot\mathrm{Sip}_{n_i}$, then
one check immediately that the subspace of $\mathcal V_{\lambda,i}$
generated by the first $\big[\frac{n_i}2\big]$ vectors of the basis $v_1^i,\ldots,v_{n_i}^i$
is $g$-isotropic. Using Remark~\ref{thm:remsignature}, we get that
$\sigma\big(g\vert_{\mathcal V_{\lambda,i}}\big)=0$ if $n_i$ is even, and
that $\vert \sigma\big(g\vert_{\mathcal V_{\lambda,i}}\big)\vert=1$ if
$n_i$ is odd.

\noindent\ \
The last statement concerning the nondegeneracy of $g\vert_{\Ker(A-\lambda)}$ follows
immedia\-tely from part (c) of Proposition~\ref{thm:restrgeigen}.
\end{proof}
\subsection{Jordan signatures}\label{sub:Jordansignature}
We will now introduce the notion of {\em Jordan signatures}, which are nonnegative integer
invariants associated to a triple $(g,A,\lambda)$, where $g$ is a nondegenerate
symmetric bilinear form on $\R^n$, $A$ is a $g$-symmetric endomorphism of $\R^n$
and $\lambda$ is an eigenvalue of $A$. For the purposes of this paper, we
will consider only the case that $\lambda$ is real.
Given such a triple $(g,A,\lambda)$, write $\Ker(A-\lambda)^n=\bigoplus_{i=1}^r\mathcal V_{\lambda,i}$
as in Proposition~\ref{thm:restrgeigen}, set $n_i=\Dim\left(\mathcal V_{\lambda,i}\right)$, denote by
$v_1^i,\ldots,v^i_{n_i}$ a basis of $\mathcal V_{\lambda,i}$ as in part (c) of
Proposition~\ref{thm:restrgeigen}
and let $\epsilon_i\cdot\mathrm{Sip}_{n_i}$ be the matrix representation of $g\vert_{\mathcal V_{\lambda,i}}$
relatively to this basis.
For $i=1,\ldots,r$, define $\varsigma_i(g,A,\lambda)$ to be the index
of the restriction of $g$ to $\mathcal V_{\lambda,i}$, and define $\varrho_i(g,A,\lambda)$ as the index of the (degenerate)
symmetric  bilinear form $b_{\lambda,i}$ on $\mathcal V_{\lambda,i}$ whose matrix
representation in the given basis is:
\begin{equation}\label{eq:defblambdai}
b_{\lambda,i}\cong\begin{pmatrix}0&0&0&\ldots&0&0&0\\ 0&0&0&\ldots&0&0&\epsilon_i\\
0&0&0&\ldots&0&\epsilon_i&0\\0&0&0&\ldots&\epsilon_i&0&0\\ &&&\vdots\\ 0&0&\epsilon_i&\ldots&0&0&0\\ 0&\epsilon_i&0&\ldots&0&0&0
\end{pmatrix}.
\end{equation}
Finally, set $\tau_i(g,A,\lambda)=\varrho_i(g,A,\lambda)+1-\varsigma_i(g,A,\lambda)$.
With the help of Remark~\ref{thm:remsignature}, such numbers can be computed explicitly as follows:
\begin{equation}\label{eq:defvarsigmai}
\varsigma_i(g,A,\lambda)=\begin{cases}\dfrac{n_i}2,\qquad&\text{if $n_i$ is even;}\\[10pt]
\dfrac{n_i-1}2,\qquad&\text{if $n_i$ is odd and $\epsilon_i>0$;}\\[10pt]
\dfrac{n_i+1}2,\qquad&\text{if $n_i$ is odd and $\epsilon_i<0$;}
\end{cases}
\end{equation}
\begin{equation}\label{eq:defvarrhoi}
\varrho_i(g,A,\lambda)=\begin{cases}\dfrac{n_i-1}2,\qquad&\text{if $n_i$ is odd;}\\[10pt]
 \dfrac {n_i}2-1,\qquad &\text{if $n_i$ is even and $\epsilon_i>0$;}
\\[10pt]
\dfrac {n_i}2,\qquad &\text{if $n_i$ is even and $\epsilon_i<0$.}
\end{cases}
\end{equation}
and
\begin{equation}\label{eq:deftaui}
\tau_i(g,A,\lambda)=\frac{1+\epsilon_i(-1)^{n_i+1}}2\in\{0,1\}.  
\end{equation}

\begin{defin}\label{thm:defJordansignature}
The \emph{Jordan signatures}, $\varsigma(g,A,\lambda)$,  $\varrho(g,A,\lambda)$
and $\tau(g,A,\lambda)$ are defined respectively as  $\sum\limits_{i=1}^r\varsigma_i(g,A,\lambda)$,
$\sum\limits_{i=1}^r\varrho_i(g,A,\lambda)$ and $\sum\limits_{i=1}^r\tau_i(g,A,\lambda)$.
\end{defin}

From \eqref{eq:deftaui}, $0\le\tau(g,A,\lambda)\le r=\Dim\big(\Ker(A-\lambda)\big)$;
moreover, $\varsigma(g,A,\lambda)$
coincides with the index of the restriction of $g$ to  $\Ker(A-\lambda)^n$,
and we get:
\begin{equation}\label{eq:reltaurhon}
\tau(g,A,\lambda)=\varrho(g,A,\lambda)+\Dim\big(\Ker(A-\lambda)\big)-\mathrm n_-\big(g\vert_{\Ker(A-\lambda)^n}\big).
\end{equation}
\end{section}

\begin{section}{Eigenvalues and conjugate points}
\label{sec:eigenvalues}
Let $I\subset\R$ be an interval and $t\mapsto \mathfrak a(t)$, $t\mapsto \mathfrak b(t)$ be continuous
maps on $I$ taking values in the space of linear endomorphisms of $\R^n$.
We can give the following general
definition:
\begin{defin}\label{thm:defconjugateinstant}
Two instants $t_0,t_1\in I$, $t_0<t_1$ are said to be {\em conjugate\/}
for the second order linear system $v''+\mathfrak a(t)v'+\mathfrak b(t)v=0$ in $\R^n$
(we also say that $t_1$ is conjugate to $t_0$) if there exists
a non identically zero solution $v$ of the system such that $v(t_0)=v(t_1)=0$.
Clearly, the set of such solutions is a
vector space whose dimension is less than or equal to $n$; such dimension is defined to be the
{\em multiplicity\/} of the conjugate instant $t_1$.
\end{defin}
\begin{rem}\label{thm:conjinstinizio}
Consider the second order linear system $v''+\mathfrak a(t)v'+\mathfrak b(t)v=0$ in $\R^n$
with $\mathfrak a,\mathfrak b:[a,b]\to\mathrm{End}(\R^n)$ continuous maps.
There exists $\varepsilon>0$ such that the set $\mathcal{C}=\big\{t\in\left]a,b\right]:t\ \text{is conjugate to}\ a\big\}$
does not contain any point of the interval $\left]a,a+\varepsilon\right]$.
To see this, consider the associated first order system in $\R^{2n}$:
$\left(\begin{smallmatrix}v\\w\end{smallmatrix}\right)'=X(t)\left(\begin{smallmatrix}v\\w\end{smallmatrix}\right)$,
with $X=\left(\begin{smallmatrix}0&\mathrm I_n\\ -\mathfrak b(t)&-\mathfrak a(t)\end{smallmatrix}\right)$,
and let $\Phi=\left(\begin{smallmatrix}\Phi_{11}&\Phi_{12}\\\Phi_{21}&\Phi_{22}\end{smallmatrix}\right):[a,b]\to\mathrm{GL}(\R^{2n})$
be its fundamental solution, i.e., $\Phi'=X\Phi$ and $\Phi(a)=\mathrm I_{2n}$.
An instant $t$ belongs to $\mathcal{C}$ iff $\Phi_{12}(t)$ is singular; since $\Phi_{12}(a)=0_n$ and
$\Phi_{12}'(a)=\mathrm I_n$, then  $\Phi_{12}(t)$ is positive definite for
$t\in\left]a,a+\varepsilon\right]$ when $\varepsilon>0$ is small enough. This proves our assertion.
\end{rem}

\begin{rem}\label{thm:remdiscrconjinstanal}
If the coefficients $\mathfrak a(t)$ and $\mathfrak b(t)$ are real analytic
functions of $t$, then it is easy to see that the set of conjugate instants
is discrete. Namely, if $v_1,\ldots,v_n$ are linearly independent
solutions of the system $v''+\mathfrak a(t)v'+\mathfrak b(t)v=0$
satisfying $v_i(t_0)=0$ for all $i$, then by standard regularity
arguments each map $v_i$ is real analytic, and the conjugate instants
correspond to the zeroes of the real analytic map $t\mapsto\mathrm{det}\big(v_1(t),\ldots,v_n(t)\big)$.
Such map is not identically zero by Remark~\ref{thm:conjinstinizio}.
\end{rem}

In case of system with constant coefficients,
the existence of conjugate instants  is related
to the spectrum of the coefficients in a quite straightforward way.
For our purposes, we will be interested in the following situation:
\begin{lem}\label{thm:existenceconjpts}
Let $A$ be an arbitrarily fixed endomorphism of $\R^n$.
There exists pairs of conjugate instants $t_0,t_1\in\R$ for
the system $v''=Av$ if and only if  $A$
has real negative eigenvalues.
\end{lem}
\begin{proof}
Since the system has constant coefficients, translations preserve its solutions,
and therefore it is not restrictive to consider the case $t_0=0$.
We consider the complexified system $v''=\complex A v$  in $\C^n$.
The first observation is that  establishing
whether an instant $t_1>0$ is conjugate to $0$ is equivalent to determining
the existence of a complex solution $v:[0,t_1]\to\C^n$ of this system which
is not identically zero and satisfying $v(0)=v(t_1)=0$. Namely, given any such
solution, its real part and its imaginary part are solutions of the real system,
and they both vanish at $0$ and at $t_1$;  at least one of the two parts cannot
vanish identically.

We can now consider a suitable basis of $\C^n$ where $\complex A$ is represented
by its Jordan form; it is immediate to see that the existence of a non
trivial solution of $v''=\complex A v$ vanishing at two given instants is equivalent
to the existence of a non trivial solution in $\C^k$ of at least one of the
systems $w''=J_k(\lambda) w$
vanishing at the same two instants. Here $\lambda$ runs in the
spectrum  of $\complex A$ and $J_k(\lambda)$ is any one of the  Jordan
blocks appearing in the Jordan decomposition of $\complex A$.

It is therefore not restrictive to assume that the spectrum of $\complex A$
consists of a single eigenvalue $\lambda\in\C$, and that $\complex A$ is represented
(in the canonical basis of $\C^n$) by the Jordan block $J_n(\lambda)$ as in \eqref{eq:Jordanblock}.
Let $v=(v_1,\ldots,v_n):\R\to\C^n$ be a non
trivial solution of $v''=J_n(\lambda)v$ vanishing at $0$ and at some other instant
$t_1>0$. Assume that the $n$-th component $v_n:\R\to\C$ of $v$ is not identically $0$;
it is easily computed $v_n=C_n\big(e^{\alpha t}-e^{-\alpha t}\big)$ for
some $C_n\in\C\setminus\{0\}$,  where $\alpha$ is any one of the two complex roots of $\lambda$.
Since $v_n(t_1)=0$, then $e^{\alpha t_1}=e^{-\alpha t_1}$, i.e., $2t_1\alpha$ is an integer
multiple of $2\pi i$, i.e., $\alpha=\frac{k\pi}{t_1} i$ for some $k\in\Z$, and therefore
$\lambda=\alpha^2$ is a negative real number. On the other hand, if $v_n$ vanishes
identically, then one computes easily $v_{n-1}=C_{n-1}\big(e^{\alpha t}-e^{-\alpha t}\big)$, to
which the same argument applies, i.e.,  $\lambda\in\R^-$ unless  $v_{n-1}$  vanishes identically.
An immediate induction argument completes
the proof: if any one of the component $v_k$ of $v$ is not identically zero,
then $\lambda\in\R^-$, and we are done. The converse is easy.
\end{proof}
By exploiting the argument in the proof of Lemma~\ref{thm:existenceconjpts}
one obtains  precise information on the displacement and
the number of conjugate instants for the system
$v''=Av$. Let us agree that by the ``number of conjugate instants'' we mean
that each conjugate instant has to be counted with its multiplicity.
\begin{cor}\label{thm:numeroconjpts}
Let $T>0$ be fixed and let $A$ be an arbitrary linear endomorphism of $\R^n$.
An instant $t_1\in\left]0,T\right]$ is conjugate to $0$ for the system $v''=Av$
if and only if there exists a real negative eigenvalue $\lambda$ of $A$ and
a positive integer $k$ such that $t_1=\frac{k\pi}{\sqrt{\vert\lambda\vert}}$.
Given such a conjugate instant $t_1$, its multiplicity is given by the sum:
\[\sum_{\lambda}
\Dim\big(\Ker(A-\lambda)\big),\]
where the sum is taken over all $\lambda$'s in the real negative spectrum of $A$
of the form $-k^2\frac{\pi^2}{t_1^2}$ for some $k\in\N\setminus\{0\}$. The number of conjugate
instants to $0$ in $\left]0,T\right]$ is given by:
\begin{equation}\label{eq:squarebrakets2}
\sum_{\lambda\in\mathfrak s(A)\cap\left]-\infty,-\frac{\pi^2}{T^2}\right]}\Dim\big(\Ker(A-\lambda)\big)\cdot
\biglip\frac{T\sqrt{\vert\lambda\vert}}{\pi}\bigrip,
\end{equation}
where $\lip\alpha\rip$ denotes the \emph{integer part} of the real number $\alpha$.
\end{cor}
\begin{proof}
Each $\lambda$-Jordan block of $A$ as in \eqref{eq:Jordanblock} gives a contribution
of $1$ to the multiplicity of the conjugate instant $t_1=\frac{k\pi}{\sqrt{\vert\lambda\vert}}$;
namely, the only non trivial solution of $v''=Av$ vanishing at $0$ and at $t_1$
when $A$ is represented by a $\lambda$-Jordan block as in \eqref{eq:Jordanblock} with $\lambda<0$ is given by
$v(t)=\big(C_1\sin(t\sqrt{\vert\lambda\vert}),0,\ldots,0\big)$,
for some $C_1\in\R$ (observe that this fact can be easily obtained from \eqref{eq:exptX} and \eqref{eq:seno}).

The conclusion follows easily from the observation that the number of
$\lambda$-Jordan blocks appearing in the Jordan form of $A$ equals
the dimension of $\Ker(A-\lambda)$.
\end{proof}

\end{section}

\begin{section}{Computation of the Maslov index}
\label{sec:computation}
Let us fix throughout this section a non degenerate
symmetric bilinear form $g$ on $\R^n$, a $g$-symmetric
linear endomorphism $A$ of $\R^n$ and a positive instant $T$; the corresponding differential
system is:
\begin{equation}\label{eq:diffsist}
v''(t)=Av(t),\qquad t\in\left[0,T\right].
\end{equation}
%
Consider the vector space $\R^n\oplus{\R^n}^*$ endowed with its canonical
symplectic form \begin{equation}\label{eq:canonicalsymplform}
\omega\big((v,\alpha),(w,\beta)\big)=\beta(v)-\alpha(w),\qquad
v,w\in\R^n,\ \alpha,\beta\in{\R^n}^*.
\end{equation}
Equation \eqref{eq:diffsist} can also be written as a first order system
in $\R^n\oplus{\R^n}^*$, using  explicitly the bilinear form $g$, as:
\begin{equation}\label{eq:forma simplectica}
\begin{pmatrix}v\\\alpha\end{pmatrix}'=\begin{pmatrix}0&g^{-1}\\
gA&0\end{pmatrix}\begin{pmatrix}v\\\alpha\end{pmatrix},
\end{equation}
from which the symplectic structure of \eqref{eq:diffsist} appears naturally (see Subsection~\ref{sub:symplecticsystem}).
The endomorphism $X=\begin{pmatrix}0&g^{-1}\\
gA&0\end{pmatrix}$ of $\R^n\oplus{\R^n}^*$ belongs to the Lie algebra of
the symplectic group $\mathrm{Sp}(\R^n\oplus{\R^n}^*,\omega)$; the fundamental solution
$\Phi(t)$ of \eqref{eq:forma simplectica} is easily computed as the exponential
$\exp(t\cdot X)$:
\begin{equation}\label{eq:exptX}
\Phi(t)=\exp(t\cdot X)=\begin{pmatrix}\mathcal C(t^2A)&t\mathcal S(t^2A)g^{-1}\\ tgA\mathcal S(t^2A)&g\,\mathcal C(t^2A)g^{-1}\end{pmatrix},
\end{equation}
where, for $B\in\mathrm{End}(\R^n)$, we have set:
\[\mathcal C(B)=\sum_{k=0}^\infty\frac1{(2k)!}B^k,\qquad \mathcal S(B)=\sum_{k=0}^\infty\frac1{(2k+1)!}B^k.\]

The conjugate instants of the system \eqref{eq:diffsist} are precisely the
instants $t$ for which the upper right block of $\Phi(t)$ is singular, i.e.,
$t\in\left]0,T\right]$ is a conjugate instant of \eqref{eq:diffsist} if and only if
$\mathcal S(t^2A)$ is singular.

Using Lemma~\ref{lem:seriediA}, we can compute explicitly $\Phi(t)$ in a Jordan basis
for $A$ using the following:
\begin{equation}\label{eq:coseno}
\mathcal C(t^2A)=\begin{pmatrix}\cos\alpha t&\frac{t\sin \alpha t}{2\alpha}&*&*&\ldots&*\\ 0&\cos\alpha t&\frac{t\sin \alpha t}{2\alpha}&*&\ldots&*\\
&&\vdots\\ 0&0&\ldots&\cos\alpha t&\frac{t\sin \alpha t}{2\alpha}&*\\ 0&0&\ldots&0&\cos\alpha t&\frac{t\sin \alpha t}{2\alpha}\cr 0&0&\ldots&0&0&\cos\alpha t\end{pmatrix}
\end{equation}
and
\begin{equation}\label{eq:seno}
\mathcal S(t^2A)=\begin{pmatrix}\frac{\sin \alpha t}{\alpha t}&\beta(t)&*&*&\ldots&*\\ 0&\frac{\sin \alpha t}{\alpha t}&\beta(t)&*&\ldots&*\\
&&\vdots\\ 0&0&\ldots&\frac{\sin \alpha t}{\alpha t}&\beta(t)&*\\ 0&0&\ldots&0&\frac{\sin \alpha t}{\alpha t}&\beta(t)\cr 0&0&\ldots&0&0&\frac{\sin \alpha t}{\alpha t}\end{pmatrix}
\end{equation}
where $\beta(t)=\frac{1}{2\alpha^2}\left(\frac{\sin \alpha t}{\alpha t}-\cos t\alpha\right)$.
\smallskip

For all $t\in\R$, the space:
\begin{equation}\label{eq:curvainLambda}
\ell(t)=\Phi(t)\big(\{0\}\oplus{\R^n}^*\big)
\end{equation}
is a Lagrangian subspace of $(\R^n\oplus{\R^n}^*,\omega)$, i.e.,
$\ell(t)$ is an $n$-dimensional subspace on which $\omega$ vanishes.
The map $t\mapsto\ell(t)$ is a real-analytic map in the Lagrangian
Grassmannian $\Lambda$ of $(\R^n\oplus{\R^n}^*,\omega)$;
given a Lagrangian $L_0$, we will denote by $\mu_{L_0}$ the
\emph{$L_0$-Maslov index}. There is a vast literature on the Maslov index, and the most
standard reference is \cite{RobSal}; we will use a slightly different
definition of Maslov index and we will follow more closely
the approach presented in \cite{GiaPicPor}. If we denote by $\Sigma_{L_0}$
the \emph{$L_0$-Maslov cycle}, which is the subset of $\Lambda$ consisting
of all Lagrangians $L$ that are not transversal to $L_0$, then
roughly speaking the $L_0$-Maslov index of a path $\ell $ is given
by the intersection number of $\ell$ and $\Sigma_{L_0}$.
When the endpoints of $\ell$ do not lie on $\Sigma_{L_0}$,
this intersection number can be computed as the
class of $\ell$ in the first relative homology group $H_1(\Lambda,\Sigma_{L_0})\cong\Z$.
The definition of the Maslov index in the general case is
as follows. Assume that $\ell:[a,b]\to\Lambda$ is a continuous
curve for which there exists $L_1\in\Lambda$ such that
$L_1\cap L_0=L_1\cap\ell(t)=\{0\}$ for all $t\in[a,b]$. Then, define
the $L_0$-Maslov index of $\ell$ as:
\begin{equation}\label{eq:primadefmul0}
\mu_{L_0}(\ell)=\mathrm n_+\big(\varphi_{L_0,L_1}(\ell(b)\big)+\Dim\big(\ell(b)\cap L_0\big)
-\mathrm n_+\big(\varphi_{L_0,L_1}\big(\ell(a)\big)-\Dim\big(\ell(a)\cap L_0\big),
\end{equation}
where, for $L\in\Lambda$ such that $L\cap L_1=\{0\}$,
$\varphi_{L_0,L_1}(L)$ is the symmetric bilinear form on $L_0$ given by
$\omega(T\cdot,\cdot)$, $T$ being the unique linear
map $T:L_0\to L_1$ whose graph \[\Gr(T)=\big\{x+Tx:x\in L_0\big\}\] is $L$.
It is not hard to prove that the right hand side of \eqref{eq:primadefmul0}
does not depend on the choice of $L_1$. Moreover, by \cite[Corollary~3.5]{GiaPicPor}, there exists
a unique extension of the $\Z$-valued map $\mu_{L_0}$ above to the set of all
continuous curves in $\Lambda$ which is invariant by fixed endpoints homotopies
and additive by concatenation.\footnote{%
We briefly observe here that our notion of Maslov index $\mu_{L_0}$ and the notion of Maslov index
$\mu_{L_0}^{\textrm{RS}}$ discussed in \cite{RobSal}, which is a half-integer, differ only in the way of counting the
contribution of the endpoints. For a continuous curve $\gamma:[a,b]\to\Lambda$,
the two quantities are related by the following simple identity:  $\mu_{L_0}^{\textrm{RS}}(\gamma)=\mu_{L_0}(\gamma)+\textstyle{1\over2}{\rm dim}\big(\gamma(a)
\cap L_0\big)-\textstyle{1\over2}{\rm dim}\big(\gamma(b)\cap L_0\big)$. In particular, if $\gamma$ has
both endpoints transversal to $L_0$,
then $\mu_{L_0}(\gamma)=
\mu_{L_0}^{\textrm{RS}}(\gamma)$.}
The Maslov index is also \emph{symplectic additive}, in the sense
that, given symplectic spaces $(V_s,\omega_s)$, Lagrangians $L_0^s\in\Lambda(V_s,\omega_s)$
and continuous paths $\ell_s:[0,T]\to\Lambda(V_s,\omega_s)$, with
$s=1,\ldots,k$, then $\mu_{\bigoplus_{s=1}^kL_0^s}(\bigoplus_{s=1}^k\ell_s)=
\sum_{s=1}^k\mu_{L_0^s}(\ell_s)$.
Finally, if $\ell:[a,b]\to\Lambda$ is a continuous path such
that $\Dim\big(\ell(t)\cap L_0\big)$ is constant on $[a,b]$, then
$\mu_{L_0}(\ell)=0.$

\begin{defin}\label{thm:defMaslovsyst}
We will denote by $\mu(g,A,T)$ the \emph{Maslov index
of the system \eqref{eq:diffsist}}, which
is defined as:
\begin{equation}\label{eq:defMaslovindexsystem}
\mu(g,A,T)=\mu_{L_0}(\ell),
\end{equation}
where $\ell:[0,T]\to\Lambda$ is the smooth curve
given in \eqref{eq:curvainLambda} and $L_0=\{0\}\oplus{\R^n}^*$.
\end{defin}
The reader should observe that, when \eqref{eq:diffsist} comes
from the Jacobi equation along a semi-Riemannian geodesic,
it is customary in the literature (see \cite{Hel}) to define the Maslov
index of the geodesic as $\mu_{L_0}\big(\ell\vert_{[\varepsilon,T]}\big)$, where
$\varepsilon>0$ is small enough so that there are no conjugate
instants of \eqref{eq:diffsist} in $\left]0,\varepsilon\right]$ (recall Remark~\ref{thm:remdiscrconjinstanal}).
The contribution to $\mu_{L_0}(\ell)$ given by the initial instant
$t=0$ is easily computed as $-\mathrm n_-(g)$ (see also Proposition~\ref{thm:contribmaslov}), so that $\mu(g,A,T)$
coincides with $\mu_{L_0}\big(\ell\vert_{[\varepsilon,T]}\big)-\mathrm n_-(g)$.
\smallskip

The intersections of the curve $\ell$ in \eqref{eq:curvainLambda}
with the $L_0$-Maslov cycle occur precisely at the conjugate
instants of the system \eqref{eq:diffsist}; when $g$ is positive definite,
then each conjugate instant gives a positive contribution to
the computation of the Maslov index, given by its multiplicity.
More generally, given a $C^1$-curve $\ell$ in $\Lambda$
which intercepts at $t=t_0$ \emph{transversally} the regular part
of the $L_0$-Maslov cycle (in which case such intersection is
isolated), the contribution to the Maslov index of $\ell$ given
by $t_0$  can be computed as the signature of a certain
symmetric bilinear form on $\ell(t)\cap L_0$ (see \cite{RobSal}).
A conjugate instant $t_0$ of the system $v''=Av$ will
be called \emph{nondegenerate} if the corresponding intersection
with $\Sigma_{L_0}$ is transverse.

The purpose of this section is to give a formula for
computing the Maslov index in the case that $g$ is arbitrary,
in which case the intersection with $\Sigma_{L_0}$ of the Lagrangian
path $\ell$ given in \eqref{eq:curvainLambda} may be degenerate
(see Corollary~\ref{thm:cordegconjinst}).

We start with the following:
\begin{lem}\label{thm:cazzatina}
Suppose that $(W_s)_{s=1}^k$ is a family of $A$-invariant and
$g$-or\-thogonal subspaces of $\R^n$ such that $\R^n=\bigoplus_{s=1}^kW_s$;
denote by $A_s:W_s\to W_s$ the restriction of $A$ to $W_s$
and by $g_s$ the restriction of $g$ to $W_s\times W_s$.

Then, $g_s$ is nondegenerate, $A_s$ is
a $g_s$-symmetric endomorphism of $W_s$ for all $s$, and
$\mu(g,A,T)=\sum_{s=1}^k\mu(g_s,A_s,T)$.
\end{lem}
\begin{proof}
It follows easily from the symplectic additivity  of the Maslov index.
Under the assumptions of the Lemma, the symplectic space $(\R^n\oplus{\R^n}^*,\omega)$
is the symplectic direct sum of the spaces $(W_s\oplus W_s^*,\omega)$,\footnote{%
Here, $W_s^*$ is identified with  $g(W_s)\subset{\R^n}^*$, i.e.,
with the subspace of ${\R^n}^*$ consisting of those linear functionals that
vanish on the $g$-orthogonal complement of $W_s$.} the
Lagrangian space $L_0$ is the direct sum of the Lagrangians $\{0\}\oplus W_s^*$,
and, by the $g$-orthogonality of the $W_s$, the curve $\ell$ is the direct sum of curves $\ell_s$ obtained
from the systems $v''=A_sv$ in $W_s$.
\end{proof}
Using Lemma~\ref{thm:sommagort}, Proposition~\ref{thm:restrgeigen} and
Lemma~\ref{thm:cazzatina}, it follows
that we may restrict our computation of the Maslov index to the
case that the spectrum of $A$ consists of a single eigenvalue
$\lambda$, which is a real negative number, that the Jordan form
of $A$ consists of a single $\lambda$-Jordan block, and that the
bilinear form $g$ is represented by a matrix of the form $\epsilon\cdot\mathrm{Sip}_n$
in the canonical basis of $\R^n$ for some $\epsilon\in\{-1,1\}$.

Restriction to this case will simplify some of the computations; the contribution
to the Maslov index given by each conjugate instant will be computed
using the following:

\begin{lem}\label{thm:formulacompMaslovind}
Let $t_1\in\left]0,T\right]$ be fixed.
If $\mathcal C(t_1^2A)$ is an isomorphism of $\R^n$, then
the Lagrangian $\ell(t_1)$ is transversal to $L_1=\R^n\oplus\{0\}$, and, for
$t$ near $t_1$,   $\varphi_{L_0,L_1}\big(\ell(t)\big)$
can be identified with the symmetric bilinear form $\mathcal B_t:\R^n\times\R^n\to\R$
given by:
\begin{equation}\label{eq:formaBt}
\mathcal B_t=t \mathcal S(t_1^2A)\mathcal C(t_1^2A)^{-1}g^{-1}.
\end{equation}
\end{lem}
\begin{proof}
Transversality of $\ell(t_1)=\Phi(t_1)\big(L_0\big)$ with $L_1$ is obviously equivalent
to the nonsingularity of the lower right block of $\Phi(t)$ (see \eqref{eq:exptX}).
Formula \eqref{eq:formaBt} is obtained by a straightforward direct calculation.
\end{proof}

Lemma~\ref{thm:formulacompMaslovind} applies if we assume that the spectrum
of $A$ consists of a single negative real number:
\begin{lem}\label{thm:isoJRn}
Assume that $\mathfrak s(A)=\{\lambda\}$, with $\lambda\in\R^-$, and
$t_1=\frac{k\pi}{\sqrt{\vert\lambda\vert}}$ for some $k\in\N$.
Then,  $\mathcal C(t_1^2A)$ is an isomorphism.
\end{lem}
\begin{proof}
Under the assumption that $\mathfrak s(A)=\{\lambda\}$, in a Jordan
basis for $A$ the $n\times n$ matrix $\mathcal C(t_1^2A)$ can be
computed explicitly as an upper triangular matrix
whose diagonal  entries are equal
to $\cos(k\pi)=(-1)^k$. Such matrix is nonsingular,
and this concludes the proof.
\end{proof}
\begin{rem}\label{thm:remdipendedallipotesi}
Observe that the conclusion of Lemma~\ref{thm:isoJRn} does not hold
in general without the assumption that the spectrum of $A$ consists
of a single eigenvalue.
\end{rem}
\begin{prop}\label{thm:contribmaslov}
Under the assumptions of Lemma~\ref{thm:isoJRn},
the contribution to the Maslov index of each conjugate instant $t_1=\frac{k\pi}{\sqrt{\vert\lambda\vert}}\in\left]0,T\right[$
of \eqref{eq:diffsist} is given by the signature of $g$.
\end{prop}
\begin{proof}
We will assume that the Jordan
form of $A$ consists of a unique $\lambda$-Jordan block, and that
$\R^n$ has a basis relative to which the matrix representation of $g$ is
of the form $\epsilon\cdot\mathrm{Sip}_n$.
All the computations that will follow are done using the matrix
representations of $A$ and $g$ in such a basis.

By Lemma~\ref{thm:formulacompMaslovind} and Lemma~\ref{thm:isoJRn}, the contribution to the
Maslov index of each conjugate instant is given by the variation of
the \emph{extended coindex} (i.e., coindex plus nullity) of the path
of symmetric bilinear forms $\mathcal B_t:\R^n\times\R^n\to\R$ given
in \eqref{eq:formaBt}. In a Jordan basis for $A$, the symmetric matrix
representing $\mathcal B_t\cong t\mathcal S(t^2A)\,\mathcal C(t^2A)^{-1}g^{-1}$ can
be computed easily  using \eqref{eq:coseno} and \eqref{eq:seno} as:
\begin{equation}\label{eq:formacalBt}
\mathcal B_t\cong\begin{pmatrix}*&*&\ldots&*&*&\psi(t)\\
 *&*&\ldots&*&\psi(t)&0\\ *&*&\ldots&\psi(t)&0&0
 \\ &&\vdots\\ *&\psi(t)&\ldots&0&0&0\\ \psi(t)&0&\ldots&0&0&0\end{pmatrix}
\end{equation}
where $\psi(t)=\frac{\epsilon}{\alpha}\tan(\alpha t)$ and $\alpha=\sqrt{\vert\lambda\vert}>0$.
If $\tan(\alpha t)>0$, then the coindex of $\mathcal B_t$ equals the coindex of $g$,
while if $\tan(\alpha t)<0$,  the coindex of $\mathcal B_t$ equals the index of $g$;
observe that $\tan(\alpha t)$ is negative (resp., positive) in a left (resp., right)
neighborhood of $t_1=\frac{k\pi}{\alpha}$.

If $t_1\in\left]0,T\right[$, then the variation of (extended) coindex of
$\mathcal B_t$ on $[t_1-\varepsilon,t_1+\varepsilon]$  is given by:
\[\mathrm n_+\left(\mathcal B_{t_1+\varepsilon}\right)-\mathrm n_+\left(\mathcal B_{t_1-\varepsilon}\right)
=\mathrm n_+(g)-\mathrm n_+(-g)=\mathrm n_+(g)-\mathrm n_-(g)=\sigma(g).\qedhere\]
\end{proof}
The formula for the jump of the extended coindex at the final instant
is a little more involved, and it requires an  analysis of
the matrix representation of the bilinear form at a conjugate instant.
Using the notations in Proposition~\ref{thm:contribmaslov},
if $t_1=\frac{k\pi}{\alpha}$ for some $k\in\N$, by direct computation involving \eqref{eq:coseno} and \eqref{eq:seno}
we get:
\begin{equation}\label{eq:Btconjinst}
\mathcal B_{t_1}\cong\begin{pmatrix}*&*&*&\ldots&*&-\frac{\epsilon k\pi}{2\alpha^3}&0\\ *&*&*&\ldots&-\frac{\epsilon k\pi}{2\alpha^3}&0&0\\
*&*&*&\ldots&0&0&0\\&&&\vdots\\ *&-\frac{\epsilon k\pi}{2\alpha^3}&0&\ldots&0&0&0\\
-\frac{\epsilon k\pi}{2\alpha^3}&0&0&\ldots&0&0&0\\0&0&0&\ldots&0&0&0
\end{pmatrix}.
\end{equation}
We are now ready for the following:
\begin{prop}\label{thm:calcolocontributofinale}
Under the assumptions of Lemma~\ref{thm:isoJRn},
if $T$ is a conjugate instant of \eqref{eq:diffsist},
i.e., if $T=\frac{k\pi}{\sqrt{\vert\lambda\vert}}$ for some $k\in\N$, then its contribution to
the Maslov index is given by the Jordan signature $\tau(g,A,\lambda)$.
\end{prop}
\begin{proof}
Using \eqref{eq:Btconjinst}, the extended coindex of $\mathcal B_{T}$ can be computed as follows. In first place,
$\Dim\left(\Ker(\mathcal B_{T}\right)=1$; moreover, we observe that the coindex of $\mathcal B_{T}$
is equal to the index of the symmetric bilinear form given in \eqref{eq:defblambdai}. Recalling
the definition of the Jordan signatures \eqref{eq:reltaurhon}, we get:
\[\begin{split}\mathrm n_+(\mathcal B_T)+&\Dim\left(\Ker(\mathcal B_T)\right)-\mathrm n_+(\mathcal B_{T-\varepsilon})\\&=
\varrho(g,A,\lambda)+\Dim\big(\Ker(A-\lambda)\big)-\mathrm n_+(-g)\\&=\varrho(g,A,\lambda)+\Dim\big(\Ker(A-\lambda)\big)-\mathrm n_-(g)=\tau(g,A,\lambda).\end{split}\]
This concludes the proof.
\end{proof}
Summarizing, we have proved the following:
\begin{cor}\label{thm:maslovindex}
For each conjugate instant $t\in\left]0,T\right]$ of \eqref{eq:diffsist}, denote
by $\mu_t(g,A)$ the contribution to the Maslov index of \eqref{eq:diffsist} given
by $t$, so that:
\[\mu(g,A,T)=\sum_{\substack{t\in\left]0,T\right]\\ t\ \text{ conjugate instant of \eqref{eq:diffsist}}}}\mu_t(g,A)-\mathrm n_-(g).\]
Then, denoting by \[\mathcal N_t=\left\{-\frac{k^2\pi^2}{t^2}:k\in\N\setminus\{0\}\right\}\subset\R^-,\]
$\mu_t(g,A)$ is computed as follows:

\[\mu_t(g,A)=\begin{cases}\sum\limits_{{\lambda\in\mathfrak s(A)\cap\mathcal N_t}}
\sigma\left(g\vert_{\Ker(A-\lambda)^{n}}\right),& \text{if}\ t<T;\\[.6truecm]
\sum\limits_{\lambda\in\mathfrak s(A)\cap\mathcal N_T}
\tau(g,A,\lambda),& \text{if}\ t=T.\qed\end{cases}\]
\end{cor}
Finally, we observe that the contribution to the Maslov index given by each conjugate instant
$t\in\left]0,T\right]$ is less than or equal to its multiplicity, due to the inequality
on the signature of $g$ proved in Corollary~\ref{thm:corsignature}, and to the inequality
on the Jordan signature $\tau$ observed at the end of Subsection~\ref{sub:Jordansignature}.

We conclude with the following observation, which relates the existence of degenerate
conjugate instants with the lack of diagonalizability of the coefficients
matrix:
\begin{cor}\label{thm:cordegconjinst}
Let $t_1\in\left]0,T\right]$ be a conjugate instant of \eqref{eq:diffsist};
then, $t_1$ is a nondegenerate conjugate instant if and only if
given any real negative eigenvalue $\lambda$ of $A$ having the form
$\lambda=-\frac{k^2\pi^2}{t_1^2}$ for some integer $k\ne0$, the algebraic
multiplicity and the geometric multiplicity of $\lambda$ coincide.
\end{cor}
\begin{proof}
Let us denote by $P_1:\R^n\oplus{\R^n}^*\to\R^n$ the projection onto
the first summand.
The conjugate instant $t_1$ is nondegenerate if and only if the restriction
of $g$ to $P_1\Phi(t)(L_0)$
is nondegenerate (see for instance \cite{GiaPicPor}), i.e., recalling
\eqref{eq:exptX}, if
and only if the restriction of $g$ is nondegenerate on the image of
$\mathcal S(t_1^2A)$. A straightforward
computations shows that such condition is equivalent to the nondegeneracy
of $g$ to $\Ker(A-\lambda)$ for each eigenvalue $\lambda$ of $A$
as in the statement of the Corollary. The conclusion follows at once
from the last statement in Corollary~\ref{thm:corsignature}.
\end{proof}
\end{section}

\begin{section}{Some examples and final remarks}
\label{sec:final}

\subsection{Conley--Zehnder   index}
The fundamental solution $t\mapsto\Phi(t)$ of a symplectic system is a smooth
curve in the symplectic group; there exists a integer invariant associated to
continuous curves in the symplectic group, which is called the {\em Conley--Zehnder\/}
index. Given $\Phi\in\mathrm{Sp}(\R^n\oplus{\R^n}^*,\omega)$, then the
graph $\mathrm{Gr}(\Phi)$ of $\Phi$ is a $2n$-dimensional subspace of $V^{4n}=(\R^n\oplus{\R^n}^*)\oplus(\R^n\oplus{\R^n}^*)$.
It is easy to see that $\mathrm{Gr}(\Phi)$ is Lagrangian relatively to the symplectic
form $\bar\omega=\omega\oplus(-\omega)$ in $V^{4n}$, where $\omega$ is as in \eqref{eq:canonicalsymplform}.
More precisely:
\begin{multline*}\bar\omega\Big[\big((v_1,\alpha_1),(v_2,\alpha_2)\big),\big((w_1,\beta_1),(w_2,\beta_2)\big)\Big]\\=
\beta_1(v_1)-\alpha_1(w_1)-\beta_2(v_2)+\alpha_2(w_2).\end{multline*}
The {\em diagonal\/} $\Delta=\big\{\big((v,\alpha),(v,\alpha)\big):v\in\R^n,\ \alpha\in{\R^n}^*\big\}\subset V^{4n}$
is also a Lagrangian space relatively to $\bar\omega$, as well as the {\em anti-diagonal\/} $\Deltao$:
\[\Deltao=\Big\{\big((v,\alpha),-(v,\alpha)\big):v\in\R^n,\ \alpha\in{\R^n}^*\Big\}.\]
\begin{defin}\label{thm:defConleyZehnder}
The \emph{Conley--Zehnder index} $\mathrm{i}_{\mathrm{CZ}}(g,A,T)$ of the system \eqref{eq:diffsist} is defined to be
the Maslov index $\mu_\Delta$ of the curve $[0,T]\ni t\mapsto\mathrm{Gr}\big(\Phi(t)\big)\in\Lambda(V^{4n},\bar\omega)$,
where $\Phi(t)=\exp(tX)$ is the fundamental solution of \eqref{eq:forma simplectica}:
\[\mathrm{i}_{\mathrm{CZ}}(g,A,T)=\mu_{\Delta}\Big([0,T]\ni t\mapsto\mathrm{Gr}\big(\Phi(t)\big)\Big).\]
\end{defin}
The Conley--Zehnder index of a symplectic system is a measure of the set
of instants $t\in[0,T]$ at which the graph of the fundamental solution $\Phi(t)$ is not transversal
to $\Delta$; observe that $\Gr\big(\Phi(t)\big)$ is transversal to $\Delta$ if and only if
$1\not\in\mathfrak s\big(\Phi(t)\big)$. The set of instants $t\in[0,T]$ at which
$\Gr\big(\Phi(t)\big)$ is not transversal to $\Delta$ may fail to be
discrete, as we state in the following lemma.
\begin{lem}\label{thm:conjinstCZ}
Assume that $\Ker(A)=\{0\}$. Then, the set of instants $t\in\left]0,T\right]$ at which $\Gr(\Phi(t))$ is not
transversal to $\Delta$ is finite, and it is given by:
\[\mathcal C=\Big\{t\in\left]0,T\right]:-\frac{4k^2\pi^2}{t^2}\in\mathfrak s(A)\ \text{for some $k\in\N\setminus\{0\}$}\Big\}.\] On the other hand, if $\,0\in\mathfrak s(A)$ such a set
coincides with the whole interval $[0,T]$.
\end{lem}
\begin{proof}
The proof follows easily from the relations $\mathfrak s(\exp(tX))=\exp (\mathfrak s(tX))$ when $t\not=0$ and
$\mathfrak s(X^2)=\mathfrak s(X)^2$ obtained from Corollary~\ref{cor:eigenFu} and
$\mathfrak s(X^2)=\mathfrak s(A)$ that comes directly.
\end{proof}
\begin{rem}\label{thm:remconjinstsubset}
Note that, in the very special case of symplectic systems of the form \eqref{eq:forma simplectica},
the set $\mathcal C$ above is a (proper) subset of the set of conjugate instants
of \eqref{eq:forma simplectica} (recall Corollary~\ref{thm:numeroconjpts}). There is in general no relation between the two sets.
\end{rem}
The Conley--Zehnder index of the fundamental solution
of a constant symplectic system is already known in the literature (see
for instance \cite[Chapter~1]{Abbo1}, computed using the rotation function
in the symplectic group. For systems of the type \eqref{eq:forma simplectica},
an alternative, direct computation can be
made using the Jordan form of $A$ and the notion of Jordan signatures.

As a consequence of the statements in the Lemma~\ref{thm:conjinstCZ}, it is convenient
to reduce the calculation to the case that $A$ is invertible.
To this aim, the following result is needed; its proof is totally
analogous to the proof of Lemma~\ref{thm:cazzatina}:
\begin{lem}\label{thm:secondacazzatina}
Under the assumptions of Lemma~\ref{thm:cazzatina}, the Conley--Zehnder index
$\mathrm{i}_{\mathrm{CZ}}(g,A,T)$ is given by the sum $\sum_s\mathrm{i}_{\mathrm{CZ}}(g_s,A_s,T)$.\qed
\end{lem}
We recall that, if $\Phi\in\mathrm{Sp}(\R^n\oplus{\R^n}^*,\omega)$ has graph
which is transversal to $\Deltao$, i.e., if $-1\not\in\mathfrak s(\Phi)$, then
if we identify $\Delta$ with $\R^n\oplus{\R^n}^*$ via the projection onto the first
coordinate,
the symmetric bilinear form $\varphi_{\Deltao,\Delta}:\R^n\oplus{\R^n}^*\times \R^n\oplus{\R^n}^*\to\R$
is given by:
\begin{equation}\label{eq:cartadoppia}
2\omega\big((\mathrm I+\Phi)^{-1}(\mathrm I-\Phi)\,\cdot\,,\,\cdot\,\big).
\end{equation}
Using the relation
${\mathcal C}(t^2A)^2={\mathrm I}+t^2A{\mathcal S}(t^2A)^2$, the matrix representation
of   (\ref{eq:cartadoppia}) is
\begin{equation}\label{eq:assomatrix}
\begin{pmatrix}
2tgA{\mathcal S}(t^2A)({\mathrm I}+{\mathcal C}(t^2A))^{-1}&0\\ 0&-2t {\mathcal S}(t^2A)({\mathrm I}+{\mathcal C}(t^2A))^{-1}g^{-1}\end{pmatrix}.
\end{equation}
\begin{lem}\label{thm:casoAinv}
Let $W\subset\R^n$ denote the $g$-orthogonal space of\ $\Ker(A^n)$,
let $\widetilde g$ denote the restriction of $g$ to $W\times W$ and let
$\widetilde A:W\to W$ denote the restriction of $A$ to $W$.
Then,
\begin{align}
\mathrm i_{\mathrm{CZ}}&(g,A,T)\notag \\ \label{eq:contributivari}&=\mathrm i_{\mathrm{CZ}}(\widetilde g,\widetilde A,T)+
\Dim\big(\Ker A\big)-\,\Dim\big(\Ker(A^n)\big)-\tau(g,A,0)\\&=\mathrm i_{\mathrm{CZ}}(\widetilde g,\widetilde A,T)
-\varrho(g,A,0)-\mathrm n_+\big(g\vert_{\Ker(A^n)}\big).\notag\end{align}
\end{lem}
\begin{proof}
By Lemma~\ref{thm:secondacazzatina}, we have
$\mathrm i_{\mathrm{CZ}}(g,A,T)=\mathrm i_{\mathrm{CZ}}(\widetilde g,\widetilde A,T)+\mathrm i_{\mathrm{CZ}}(g_0,A_0,T)$,
where $A_0:\Ker(A^n)\to \Ker(A^n)$ is the restriction of $A$ to $\Ker(A^n)$ and
$g_0$ is the restriction of $g$ to $\Ker(A^n)\times\Ker(A^n)$.
A direct computation involving Lemma~\ref{lem:seriediA} and the equality $\Ker(X)=\Ker(A)$
shows that, if we set $X_0=\left(\begin{smallmatrix}0&g_0^{-1}\\ g_0A_0&0\end{smallmatrix}\right)$
and $\Phi_0(t)=\exp(t X_0)$, then $\Dim\big(\Gr(\Phi_0(t))\cap\Delta\big)=\Dim\big(\Ker(A)\big)$
for all $t\in\left]0,T\right]$, while $\Dim\big(\Gr(\Phi_0(0))\cap\Delta\big)=\Dim(\Ker(A^n))$.
Hence $\mathrm i_{\mathrm{CZ}}(g_0,A_0,T)$ is given by
the only contribution of the initial instant $t=0$; in order to compute such contribution,
using the symplectic additivity of the Maslov index we will assume that $\Dim\big(\Ker(A)\big)=1$.
Under this assumption, the Jordan form of $A_0$ has a single $0$-Jordan block of size $k_0\times k_0$,
where $k_0=\Dim\big(\Ker(A^n)\big)$, and $g$ takes the form $g=\epsilon\cdot\mathrm{Sip}_{k_0}$, with
$\epsilon\in\{-1,1\}$;
using  \eqref{eq:coseno}, \eqref{eq:seno} and \eqref{eq:assomatrix} one computes easily
the matrix representation of the symmetric bilinear form $\varphi_{\Deltao,\Delta}\big(\Gr(\Phi_0(t))\big)$
for $t>0$ near $0$, which is the direct sum of  two $k_0\times k_0$ symmetric matrices of the form:
\[\epsilon\begin{pmatrix}0&0&0&\ldots&0&0&0\\ 0&0&0&\ldots&0&0&t\\ 0&0&0&\ldots&0&t&*\\0&0&0&\ldots&t&*&*\\&&&\vdots\\0&0&t&\ldots&*&*&*\\
0&t&*&\ldots&*&*&*
\end{pmatrix}\quad\text{and}\quad \epsilon\begin{pmatrix}*&*&*&\ldots&*&*&-t\\ *&*&*&\ldots&*&-t&0\\
*&*&*&\ldots&-t&0&0\\&&&\vdots\\ *&*&-t&\ldots&0&0&0\\ *&-t&0&\ldots&0&0&0\\-t&0&0&\ldots&0&0&0\end{pmatrix}.\]
Thus, for $t>0$ near $0$, the extended coindex of $\varphi_{\Deltao,\Delta}\big(\Gr(\Phi_0(t))\big)$ is easily
computed with the help of the Jordan signatures as:
\begin{multline*}
\Dim\big(\Gr(\Phi_0(t))\cap\Delta\big)+\mathrm n_+\left(\varphi_{\Deltao,\Delta}\big(\Gr(\Phi_0(t))\big)\right)\\=1+\mathrm n_+\big(\epsilon\cdot\mathrm{Sip}_{k_0-1}\big)+
\mathrm n_+\big(-\epsilon\cdot\mathrm{Sip}_{k_0}\big)=k_0-\rho(g,A,0)+\varsigma(g,A,0).\end{multline*}
For $t=0$, the extended coindex of $\varphi_{\Deltao,\Delta}\big(\Gr(\Phi_0(0))\big)=\varphi_{\Deltao,\Delta}(\Delta)=\{0\}$ is equal to $2k_0$.
Formula \eqref{eq:contributivari} follows readily using \eqref{eq:reltaurhon}.
\end{proof}
Lemma~\ref{thm:casoAinv} tells us that, in order to compute the Conley--Zehnder index
of \eqref{eq:diffsist}, it suffices to consider the case that $A$ is invertible.

We are now ready for the following:
\begin{prop}\label{thm:compConZeh}
The contribution to the Conley--Zehnder index of \eqref{eq:diffsist} given by the
initial instant $t=0$ is given by the following formula:
\begin{multline}\label{eq:contrinizCZ}
-2\mathrm n_+(g)+\Dim\big(\Ker A\big)+\,\sigma\big(g\vert_{\Ker(A^n)}\big)-\tau(g,A,0)\\=
-2\mathrm n_+(g)-\varrho(g,A,0)+\mathrm n_+\big(g\vert_{\Ker(A^n)}\big).
\end{multline}
If $t_1\in\left]0,T\right[\cap \mathcal C$, then its contribution
to the Conley--Zehnder index of \eqref{eq:diffsist} is given by:
\[-2\sum_{\lambda}\sigma\big(g\vert_{\Ker(A-\lambda)^n}\big),\]
where the sum is taken over all $\lambda\in\mathfrak s(A)\cap\R^-$ of the form
$\lambda=-\frac{4k^2\pi^2}{t_1^2}$ for some $k\in\N\setminus\{0\}$. The contribution of
the final instant $T$ is given by
\begin{equation}
2\sum_\lambda\big(-\varrho(g,A,\lambda)+\mathrm n_-(g\vert_{\Ker(A-\lambda)^n})\big)
=2\sum_\lambda\big(-\tau(g,A,\lambda)+\dim(\Ker(A-\lambda))\big)
\end{equation}
where the sum is taken over all $\lambda\in\mathfrak s(A)\cap\R^-$ of the form
$\lambda=-\frac{4k^2\pi^2}{T^2}$ for some $k\in\N\setminus\{0\}$.

\end{prop}
\begin{proof}
The contribution to the Conley--Zehnder index given by the initial instant
of the null eigenvalue of $A$ is computed in Lemma~\ref{thm:casoAinv}.
We need to compute the contribution to the index given by the initial
instant of the \emph{reduced} symplectic system, i.e., the system  in $W\oplus W^*$ with coefficient matrix:
\[\widetilde X=\begin{pmatrix}0&\tilde g^{-1}\\\tilde g\tilde A&0\end{pmatrix},\]
where $W$ is the $g$-orthogonal subspace
to $\Ker(A^n)$ and $\tilde g$ is the restriction of $g$ to  $W$.
By Lemma \ref{thm:secondacazzatina}, we can assume that $\mathfrak s(A)=\{\lambda\}$ and
the Jordan form of $A$ consists of a single block, i.e., $\Dim\big(\Ker(A-\lambda)\big)=1$.
If $r_0=\dim(\Ker(A-\lambda)^n)$, the matrix representation of the symmetric bilinear form $\varphi_{\Deltao,\Delta}\big(\Gr(\Phi(t))\big)$
when  $\sin t\alpha\not=0$ is the direct sum of  two $r_0\times r_0$ symmetric matrices that
can be computed using again equations  \eqref{eq:coseno}, \eqref{eq:seno} and \eqref{eq:assomatrix} as:
\begin{equation}\label{eq:triangular1}\epsilon\begin{pmatrix}0&0&0&\ldots&0&0&\mu_1\\ 0&0&0&\ldots&0&\mu_1&*\\ 0&0&0&\ldots&\mu_1&*&*\\0&0&0&\ldots&*&*&*\\&&&\vdots\\0&\mu_1&*&\ldots&*&*&*\\
\mu_1&*&*&\ldots&*&*&*
\end{pmatrix}\end{equation}
and
\begin{equation}\label{eq:triangular2}
 \epsilon\begin{pmatrix}*&*&*&\ldots&*&*&\mu_2\\ *&*&*&\ldots&*&\mu_2&0\\
*&*&*&\ldots&\mu_2&0&0\\&&&\vdots\\ *&*&\mu_2&\ldots&0&0&0\\ *&\mu_2&0&\ldots&0&0&0\\\mu_2&0&0&\ldots&0&0&0\end{pmatrix},
\end{equation}
where $\mu_1=-\frac{2 \alpha\sin t\alpha}{1+\cos t\alpha}$ and
$\mu_2=-\frac{2\sin t\alpha}{\alpha(1+\cos t\alpha)}$.
From \eqref{eq:triangular1} and \eqref{eq:triangular2} it is easy to see that the contribution of the initial
instant is
\[2 \mathrm n_-(\tilde g)-2r_0=-2 \mathrm n_+(\tilde g)=2(\mathrm n_+\big(g\vert_{\Ker (A^n)}\big)-\mathrm n_+(g))\]
and formula \eqref{eq:contrinizCZ} follows now easily from \eqref{eq:contributivari}.
\smallskip

Assume now  $t_1\in\left]0,T\right[\cap\mathcal C$ and
$\lambda=-\frac{4k^2\pi^2}{t_1^2}$ for some $k\in\N\setminus\{0\}$. Then as the matrix representation
of $\varphi_{\Deltao,\Delta}\big(\Gr(\Phi(t))\big)$ is the direct sum of \eqref{eq:triangular1} and \eqref{eq:triangular2} it is easy to see that the contribution of $t_1$ is $-2\sigma (\tilde g\vert_{\Ker(A-\lambda)^n})$.

In order to compute the contribution of the final instant first we observe that the
matrix representation of $\varphi_{\Deltao,\Delta}\big(\Gr(\Phi(T))\big)$ when $T=\frac{2k\pi}{\alpha}$
is the direct sum of
the two matrices:
\[\epsilon\begin{pmatrix}0&0&0&\ldots&0&0&0\\ 0&0&0&\ldots&0&0&\frac{k\pi}{\alpha}\\ 0&0&0&\ldots&0&\frac{k\pi}{\alpha}&*\\0&0&0&\ldots&\frac{k\pi}{\alpha}&*&*\\&&&\vdots\\0&0&\frac{k\pi}{\alpha}&\ldots&*&*&*\\
0&\frac{k\pi}{\alpha}&*&\ldots&*&*&*
\end{pmatrix},\qquad
\epsilon\begin{pmatrix}*&*&*&\ldots&*&\frac{k\pi}{\alpha^3}&0\\ *&*&*&\ldots&\frac{k\pi}{\alpha^3}&0&0\\
*&*&*&\ldots&0&0&0\\&&&\vdots\\ *&\frac{k\pi}{\alpha^3}&0&0&0&\ldots&0\\ \frac{k\pi}{\alpha^3}&0&0&0&0&\ldots&0\\0&0&0&\ldots&0&0&0\end{pmatrix}.\]
Then the contribution of the final instant is
\begin{align*}
&\mathrm
n_+\left(\varphi_{\Deltao,\Delta}\big(\Gr(\Phi(T))\big)\right)+\\
&\qquad\qquad
\dim\left(\Ker(\varphi_{\Deltao,\Delta}\big(\Gr(\Phi(T))\big))\right)\!-\!
\mathrm n_+\left(\varphi_{\Deltao,\Delta}\big(\Gr(\Phi(T-\varepsilon))\big)\right)\\
&\qquad=2(\dim(\Ker(A-\lambda)^n)-\varrho(g,A,\lambda)-\mathrm n_+(g\vert_{\Ker(A-\lambda)^n}\big)\\
&\qquad=2(-\varrho(g,A,\lambda)+\mathrm n_-(g\vert_{\Ker(A-\lambda)^n}\big)\\
&\qquad=2(-\tau(g,A,\lambda)+\dim(\Ker(A-\lambda)))
\end{align*}
This concludes the proof.
\end{proof}

\subsection{Maslov index of an arbitrary constant symplectic system}
\label{sub:symplecticsystem} A more general class of differential
systems where the notion of Maslov index is naturally defined
consists of the so called symplectic systems. Denote by
$\mathrm{Sp}(2n,\R)$ the Lie group consisting of all isomorphisms
of $\R^{n}\oplus{\R^n}^*$ preserving the canonical symplectic form
of $\R^{n}\oplus{\R^n}^*$, and let $\mathrm{sp}(2n,\R)$ be its Lie
algebra. Recall that a $(2n)\times(2n)$ real matrix $X$ belongs to
$\mathrm{sp}(2n,\R)$ if and only if $X$ is written in $n\times n$
blocks as $X=\left(\begin{smallmatrix} \mathfrak A&\mathfrak
B\\\mathfrak C&-\mathfrak A^*\end{smallmatrix}\right)$, where
$\mathfrak B$ and $\mathfrak C$ are symmetric matrices. We call a
\emph{symplectic differential system} a first order system in
$\R^{n}\oplus{\R^n}^*$ of the form:
\begin{equation}\label{eq:symplsyst}
\begin{pmatrix}v\\ \alpha\end{pmatrix}'=X\begin{pmatrix}v\\\alpha\end{pmatrix}
\end{equation}
where $X:[a,b]\to\mathrm{sp}(2n,\R)$ is a smooth curve. The fundamental solution
of a symplectic differential system is a smooth curve $\Phi$ taking values in
$\mathrm{Sp}(2n,\R)$; the \emph{Maslov index of a symplectic system}
is defined to be the $L_0$-Maslov index of the curve $t\mapsto\Phi(t)(L_0)\in\Lambda$,
where $L_0$ is the Lagrangian subspace $\{0\}\oplus{\R^n}^*$ of $\R^{n}\oplus{\R^n}^*$.
Similarly, an instant $t_0\in\left]a,b\right]$ is defined to be \emph{conjugate}
for the system~\eqref{eq:symplsyst} if $\Phi(t_0)(L_0)$ belongs to the $L_0$-Maslov cycle
$\Sigma_{L_0}$; equivalently, $t_0$ is conjugate if there exists a non trivial solution
$\begin{pmatrix}v\\ \alpha\end{pmatrix}$ of \eqref{eq:symplsyst} such that $v(a)=v(t_0)=0$.
For example, the second order system \eqref{eq:diffsist} in $\R^n$ is equivalent
to the symplectic system \eqref{eq:forma simplectica} in $\R^n\oplus{\R^n}^*$ whose coefficient
matrix is the constant curve $X=\left(\begin{smallmatrix}0&g^{-1}\\gA&0\end{smallmatrix}\right)$;
the notion of Maslov index and conjugate instant for such symplectic system obviously coincide with the
corresponding notions for the system \eqref{eq:diffsist} given in
Section~\ref{sec:computation}.

We will  now show how to reduce the computation of the Maslov index of a
class of constant symplectic systems to the case of a system of the form
\eqref{eq:forma simplectica}, for which the theory developed earlier applies.
To this aim, let us fix an element $X=\left(\begin{smallmatrix}
\mathfrak A&\mathfrak B\\\mathfrak C&-\mathfrak A^*\end{smallmatrix}\right)\in\mathrm{sp}(2n,\R)$
and let us consider the corresponding symplectic system as in \eqref{eq:symplsyst}.
We want to restrict our attention to those symplectic systems
for which the set of conjugate instants is discrete, and for this
we need the following:
\begin{lem}\label{thm:conjinstdiscr}
Consider a constant symplectic system with matrix coefficients $X=\left(\begin{smallmatrix}
\mathfrak A&\mathfrak B\\\mathfrak C&-\mathfrak A^*\end{smallmatrix}\right)$ on the
interval $[0,T]$.
If the upper right $n\times n$ block $\mathfrak B$ of $X$ is non singular, the set of conjugate instants
is finite.

Conversely, if $\Ker(\mathfrak A^*)\cap\Ker(\mathfrak B)\ne\{0\}$, then every $t\in\left]0,T\right]$ is conjugate.
\end{lem}
\begin{proof}
If $\mathfrak B$ is non singular, then the conjugate instants
of the symplectic system correspond to the conjugate instants
of the second order equation in $\R^n$:
\begin{equation}\label{eq:abovequation}
v''+(\mathfrak B\mathfrak A^*\mathfrak B^{-1}-\mathfrak A)v'-(\mathfrak B\mathfrak C+\mathfrak B\mathfrak A^*\mathfrak B^{-1}
\mathfrak A)v=0.
\end{equation}
Namely, given a solution $v$ of \eqref{eq:abovequation}, the pair $(v,\alpha)$ with
 $\alpha=\mathfrak B^{-1}(v'-\mathfrak Av)$ is a solution of the symplectic system, and
 this gives a bijective correspondence between the solutions of the first order system
 and the solutions of \eqref{eq:abovequation}.
 As observed in Remark~\ref{thm:remdiscrconjinstanal}, the conjugate instants of
\eqref{eq:abovequation} form a discrete set.
Conversely, if $\alpha_0\in\Ker(\mathfrak A^*)\cap\Ker(\mathfrak B)$ is non zero,
then the constant $\left(\begin{smallmatrix}v\\ \alpha\end{smallmatrix}\right)
\equiv \left(\begin{smallmatrix}0\\ \alpha_0\end{smallmatrix}\right)$ is a non trivial
solution of the system for which $v(a)=v(t_0)=0$ for all $t_0\in\left]0,T\right]$.
\end{proof}

In view of the result above, let us now restrict our attention to those
constant symplectic systems whose coefficient matrix $X=\left(\begin{smallmatrix}
\mathfrak A&\mathfrak B\\\mathfrak C&-\mathfrak A^*\end{smallmatrix}\right)$ has
non singular upper right block $\mathfrak B$. Let us also assume
that the linear map $\mathfrak B^{-1}\mathfrak A:\R^n\to{\R^n}^*$ is
self-adjoint, i.e., that $\mathfrak B^{-1}\mathfrak A=\mathfrak A^*\mathfrak B^{-1}$,
and let us denote by $\phi$ the endomorphism of $\R^n\oplus{\R^n}^*$ which is
written in $n\times n$ blocks as:
\[\phi=\begin{pmatrix} \mathrm I&0\\\mathfrak B^{-1}\mathfrak A&\mathrm I\end{pmatrix}.\]
An immediate computation shows that, since $\mathfrak B^{-1}\mathfrak A$ is symmetric,
$\phi\in\mathrm{Sp}(2n,\R)$; moreover, $\phi(L_0)=L_0$. We compute:
\[\begin{split}\widetilde X=\phi X\phi^{-1}\;&=\begin{pmatrix} \mathrm I&0\\\mathfrak B^{-1}\mathfrak A&\mathrm I\end{pmatrix}
\begin{pmatrix}
\mathfrak A&\mathfrak B\\\mathfrak C&-\mathfrak A^*\end{pmatrix}
\begin{pmatrix} \mathrm I&0\\-\mathfrak B^{-1}\mathfrak A&\mathrm I\end{pmatrix}\\ &=\begin{pmatrix}
0&\mathfrak B\\ \mathfrak B^{-1}\mathfrak A^2+\mathfrak C&0\end{pmatrix}=\begin{pmatrix}
0&\mathfrak B\\ \mathfrak A^*\mathfrak B^{-1}\mathfrak A+\mathfrak C&0\end{pmatrix}\in\mathrm{sp}(2n,\R).
\end{split}\]
Moreover, if $\Phi(t)$ is the fundamental solution of the symplectic system
with constant coefficient matrix $X$, the fundamental solution of the symplectic
system with coefficient matrix $\widetilde X$ is easily computed as:
\[\widetilde\Phi(t)=\phi\Phi(t)\phi^{-1}.\]
Since $\phi$ preserves the Lagrangian $L_0$, the conjugate instants
of the symplectic systems corresponding to the coefficient matrices $X$
and $\widetilde X$ coincide; moreover, the symplectic invariance of
the Maslov index implies that also the Maslov indices of the two
systems coincide. Let us denote by $A$ the endomorphism of $\R^n$ given
by:
\begin{equation}\label{eq:AfrakABC} A=\mathfrak B\big(\mathfrak A^*\mathfrak B^{-1}\mathfrak A+
\mathfrak C\big)\end{equation}
and by $g$ the nondegenerate symmetric bilinear form on $\R^n$ given by:
\begin{equation}\label{eq:gfrakB}
g=\mathfrak B^{-1}.
\end{equation}

We have proven the following:
\begin{cor}\label{thm:conjinstMaslovindexsymplsystem}
Let $X=\left(\begin{smallmatrix}
\mathfrak A&\mathfrak B\\\mathfrak C&-\mathfrak A^*\end{smallmatrix}\right)\in\mathrm{sp}(2n,\R)$
be such that $\mathfrak B$ is non singular and such that $\mathfrak B^{-1}\mathfrak A$ is
symmetric. Then, denoting by $A$ the endomorphism of $\R^n$ given in
\eqref{eq:AfrakABC} and by $g$ the nondegenerate symmetric bilinear form
on $\R^n$ given by \eqref{eq:gfrakB}, the conjugate instants and
the Maslov index of the symplectic system \eqref{eq:symplsyst}
coincide respectively with the conjugate instants and the Maslov index
of the second order differential system $v''=Av$, computed in
Corollary~\ref{thm:numeroconjpts} and in Corollary~\ref{thm:maslovindex}.\qed
\end{cor}
A similar reduction is clearly possible for the computation of the
Conley--Zehn\-der index of an arbitrary constant symplectic system.
\begin{subsection}{Bi-invariant metrics on Lie groups}
As a special case of semi-Riemannian locally symmetric manifold, in
this section we will consider the case of a Lie group $G$ endowed with
a bi-invariant semi-Riemannian metric $h$. We will denote by $\mathfrak g$ the
Lie algebra of $G$; recall that a (nondegenerate) symmetric bilinear
form $h$ on $\mathfrak g$ is bi-invariant if and only if $h\big(\mathrm{ad}_XY,Z\big)=
-h(Y,\mathrm{ad}_XZ\big)$ for all $X,Y,Z\in\mathfrak g$, where $\mathrm{ad}_XY=[X,Y]$.
A description
of Lie algebras admitting semi-Riemannian bi-invariant metrics   can be found
in \cite{medina}.
\smallskip

Let us start with the following technical result:
\begin{lem}\label{thm:lemmatecnico}
Let $\mathfrak g$ be a real $n$-dimensional Lie algebra endowed with a bi-invariant nondegenerate
symmetric bilinear form $h$. Let $X_1,\ldots,X_n$ be an
$h$-orthonormal basis of $\mathfrak g$, and set:
\[[X_i,X_j]=\sum_kC^k_{ij}X_k,\ \epsilon_i=h(X_i,X_i)\in\{\pm1\},\ \epsilon=\epsilon_1\cdot\ldots\cdot\epsilon_n,\ \text{and}\  a_{ij}=
\epsilon_iC_{nj}^i,\]
for all $i,j=1,\ldots,n$. Assume that for all choice of (pairwise distinct) indices
$i,j,k,l\in\{1,\ldots,n-1\}$
the following identities hold:
\begin{equation}\label{eq:identities}
C_{ij}^nC_{kl}^n+C_{jk}^nC_{il}^n+C_{ki}^nC_{jl}^n=0.
\end{equation}
Then, the characteristic polynomial
$P(\lambda)$ of the
linear operator $\mathrm{ad}_{X_n}:\mathfrak g\to\mathfrak g$  is given by:
\begin{equation}\label{eq:polcarat}
P(\lambda)=(-1)^{n-1}\lambda^{n-3}\Big(\lambda^2+\sum_{i<j}\epsilon_i\epsilon_j a_{ij}^2\Big).\end{equation}
In the above situation, if $\alpha^2=\sum\limits_{i<j}\epsilon_i\epsilon_j a_{ij}^2>0$, then, denoting
by $E_{\pm i\alpha}$ the (complex) eigenspace of $\mathrm{ad}_{X_n}$ corresponding
to the eigenvalue $\pm i\alpha$ and by $W_\alpha\subset\mathfrak g$ the real part of $E_{i\alpha}\oplus
E_{-i\alpha}$ (which is a $2$-dimensional subspace of $\mathfrak g$), the restriction
$h\vert_{W_\alpha\times W_\alpha}$ is either positive or negative
definite.
\end{lem}
\begin{rem}
We observe that if any two of the indices $i,j,k,l$ are equal, then
the identities \eqref{eq:identities} hold automatically; this is easily checked
using the anti-symmetry properties satisfied by the coefficients $C_{ij}^k$.
From this observation, it follows immediately that the technical assumption
\eqref{eq:identities} is satisfied when $n=\mathrm{dim}(\mathfrak g)\le 4$. Moreover,
the Jacobi identity satisfied by $[\cdot,\cdot]$ is equivalent to:
\begin{equation}\label{eq:idJacobi}
\sum_m\epsilon_m\left(C^m_{ij}C^m_{kl}+C_{jk}^mC_{il}^m+C^m_{ki}C^m_{jl}\right)=0.
\end{equation}
From this equality it follows that also in the case that $\mathrm{dim}(\mathfrak g)=5$,
the assumption \eqref{eq:identities} is satisfied. When $n\ge6$, the
identities \eqref{eq:identities}
are not necessarily satisfied; for instance, they are not satisfied in the
case of the product $S^3\times S^3$ endowed with the semi-Riemannian
bi-invariant metric $h=h_0\oplus(-h_0)$, where $h_0$ is the round metric
on $S^3$.
Finally, we observe that if the metric $h$ is Lorentzian, i.e., if $h$ has index
$1$, then in the
last statement of Lemma~\ref{thm:lemmatecnico} we can in fact conclude that the
restriction $h\vert_{W_\alpha\times W_\alpha}$ is always positive definite.
\end{rem}
\begin{proof}[Proof of Lemma~\ref{thm:lemmatecnico}]
In the basis $X_1,\ldots,X_n$, the matrix representing $\mathrm{ad}_{X_n}$ is
given by $B=(C^i_{nj})_{i,j=1}^n$; since the last column and the last row of this matrix
are zero, we will  consider the  square  matrix of order $n-1$ obtained removing the last row and the last column.
In order to compute the characteristic polynomial of $B$, we observe that if
we set $A=(a_{ij})$ then
\[{\det}(B-\lambda I)=\epsilon\epsilon_n\det (A-D)\]
where $D$
is a diagonal matrix of order $n-1$ with diagonal elements $d_{ii}=\epsilon_i \lambda$. Since $A$
is skew-symmetric, its determinant is zero when $n-1$ is odd; if $n-1$ is even
then the determinant of $A$ can be computed as the square of
the {\em Pfaffian\/} of $A$.
Let us recall briefly the notion of Pfaffian. Let $\pi=\{(i_1,j_1),(i_2,j_2),\ldots,(i_r,j_r)\}$
be a partition
of $\{1,\ldots,n-1\}$ with $i_k<j_k$ for $k=1,\ldots,r=\frac{n-1}2$, where the order of the pairs
is not taken into account.
We set
$a_\pi=a_{i_1j_1}a_{i_2j_2}\ldots a_{i_rj_r}$ and we denote by $\alpha_\pi$
the permutation
$(i_1j_1i_2j_2\ldots i_rj_r)$; the Pfaffian of $A$ is then defined
as
\[\sum_\pi {\rm sg}(\alpha_\pi)a_\pi.\]
In order to get an expression for the
characteristic polynomial of $B$, we define $\{k_1\ldots k_{2t}\}$ with $k_1<k_2<\ldots <k_{2t}$,
as the Pfaffian of
the matrix obtained by taking the rows and the columns $k_1,\ldots,k_{2t}$ of the
matrix $A$. Then
\begin{align*}
P(\lambda)&={\rm Det}(B-\lambda I)=\epsilon\epsilon_n {\rm Det}(A-D)\\
&=\sum_{t=0}^{[\frac{n-1}{2}]}\Big(\sum_{
k_1<\ldots<k_{2t}}\{k_1\ldots k_{2t}\}^2\epsilon_{k_1}\ldots\epsilon_{k_{2t}}\Big)(-\lambda)^{n-2t-1}
\end{align*}
where $k_1,\ldots,k_{2t}$ run in the set $\{1,\ldots,n-1\}$.

Now,  the identity (\ref{eq:identities}) is equivalent to $\{k_1k_2k_3k_4\}=0$ for all
$k_1,k_2,k_3,k_4\in \{1,\ldots,n-1\}$, hence we get:
\begin{equation}\label{eq:dadimostrare}
P(\lambda)=(-\lambda)^{n-3}\Big(\lambda^2+\sum_{i<j}\epsilon_i\epsilon_j a_{ij}^2\Big).
\end{equation}
To show this we observe that the following relation holds:
\begin{multline*}
\{k_1\ldots k_{2t}\}=\sum_{m<n}{\rm sg}(\alpha_{m,n})\{k_1k_2mn\}
\{[k_3k_4\ldots k_{2t}]_{m,n}\}\\-(t-2)a_{k_1k_2}\{k_3k_4\ldots k_{2t}\}
\end{multline*}
where $m,n$ take values in $\{k_3,\ldots,k_{2t}\}$,
$[k_3k_4\ldots k_{2t}]_{m,n}$ denotes the ordered set obtained by removing
$m$ and $n$ from the list $(k_3,k_4,\ldots, k_{2t})$, and $\alpha_{m,n}$ is the permutation
\[(k_1k_2mn[k_3k_4\ldots k_{2t}]_{m,n}).\] Formula~\eqref{eq:dadimostrare} is obtained now
using induction and
the identities (\ref{eq:identities}).

When
$\sum_{i<j}\epsilon_i\epsilon_j a_{ij}^2\not=0$, at least one of the coefficients
$a_{ij}$ is non null and there exists $p\in\{1,\dots,n-1\}$ such that $\sum_i \epsilon_ia_{i p}^2\not=0$;
for simplicity we will assume that  $a_{12}$ is not zero.
In this case, by using the identities \eqref{eq:identities} is easy to see that the system
\[T_j=(a_{2j},-a_{1j},0,\ldots,0,
a_{12},0,\ldots,0),\quad j=3,\ldots,n-1,\] where $a_{12}$ appears in the
$j$-th position, is a basis  of the Eigenspace associated to the zero Eigenvalue of $A$.
Using again
the identities (\ref{eq:identities}) we get that the vectors
\begin{align*}
P&=(\epsilon_1 a_{1 p},\epsilon_2 a_{2 p},
\ldots,\epsilon_i a_{i p},\ldots,\epsilon_{n-1}a_{n-1 p}),\\
Q&=(\epsilon_1 \sum_i\epsilon_i a_{1i}a_{i p},
\epsilon_2 \sum_i\epsilon_i a_{2i}a_{i p},
\ldots,\epsilon_{n-1}\sum_i\epsilon_i a_{n-1 i}a_{i p}),
\end{align*}
are orthogonal to every $T_j$. Furthermore
\begin{align*}
g(P,P)&=\sum_i \epsilon_ia_{i p}^2\not=0,\\
g(P,Q)&=0,\\
g(Q,Q)&=\Big(\sum_i\epsilon_ia_{ip}^2\Big)\Big(\sum_{i<j}\epsilon_i\epsilon_j
a_{ij}^2\Big)\not=0.
\end{align*}
Therefore, the system $\{P,Q\}$ is a basis of $W_\alpha$.
Then, the signature of the
restriction of $h$ to $W_\alpha$ is computed looking at the sign of the expression:
\begin{equation*}g(P,P)g(Q,Q)-g(P,Q)^2=
\big(\sum_i\epsilon_ia_{ip}^2\big)^2\big(\sum_{i<j}\epsilon_i\epsilon_j
a_{ij}^2\big),
\end{equation*}
which is positive, from which the last statement of the lemma follows.
\end{proof}

The geodesics through the identity of a Lie group $G$ endowed with
a bi-invariant semi-Riemannian metric $h$ are the one-parameter subgroups
of $G$. The covariant derivative of the Levi--Civita connection is given, in the
case of left-invariant vector fields $X,Y$ on $G$, by:
\begin{equation}\label{eq:curvatura}
\nabla_XY=\tfrac12\mathrm{ad}_XY=\tfrac12[X,Y],
\end{equation}
and, for $X,Y,Z\in\mathfrak g=T_1G$, the curvature tensor is given by:
\[R_{XY}Z=\nabla_X\nabla_YZ-\nabla_Y\nabla_XZ-\nabla_{[X,Y]}Z=\tfrac14[Z,[X,Y]]=\tfrac14\mathrm{ad}_Z\mathrm{ad}_XY.\]
\begin{prop}\label{thm:gruppi}
Let $G$ be an $n$-dimensional real Lie group endowed with a bi-invariant semi-Riemannian
metric tensor $h$ and let $\mathfrak g$ be its Lie algebra. Let $\gamma:\R\to G$ be a one-parameter subgroup
of $G$ with $X=\gamma'(0)$,  and let $t_0\in\left]0,+\infty\right[$ be fixed.
Then:
\begin{itemize}
\item[(a)] $\gamma(t_0)$ is conjugate to $\gamma(0)=1$ along $\gamma$ if and only if the spectrum
$\mathfrak s\left(\mathrm{ad}_X\right)$
of the linear operator $\mathrm{ad}_{X}:\mathfrak g\to\mathfrak g$ contains
a purely imaginary number of the form $2ki\pi t_0^{-1}$ for some $k\in\N\setminus\{0\}$.
\end{itemize}
If $\gamma(t_0)$ is conjugate to $\gamma(0)$ along $\gamma$, set $\mathcal K_{t_0}=\Big\{k\in\N:
2ki\pi t_0^{-1}\in\mathfrak s\left(\mathrm{ad}_X\right)\Big\}$,
\[W_{t_0}=\bigoplus_{k\in\mathcal K_{t_0}}\Ker\Big(\mathrm{ad}_X^2+\frac{4k^2\pi^2}{t_0^2}\Big),
\quad \widetilde W_{t_0}=\bigoplus_{k\in\mathcal K_{t_0}}\Ker\Big(\mathrm{ad}^2_X+\frac{4k^2\pi^2}{t_0^2}\Big)^n. \]
Then:
\begin{itemize}
\item[(b)] $\gamma(t_0)$ is nondegenerate if and only if $W_{t_0}=\widetilde W_{t_0}$;
\smallskip

\item[(c)] the multiplicity of $\gamma(t_0)$ is given by $\mathrm{dim}\left(W_{t_0}\right)$, which is an even
number;
\smallskip

\item[(d)]  the contribution of $\gamma(t_0)$ to the Maslov index is
$\sigma\Big(h\vert_{\widetilde W_{t_0}\times\widetilde W_{t_0}}\Big)$;
\smallskip

\item[(e)] if $\sigma\left(h\vert_{\widetilde W_{t_0}\times\widetilde W_{t_0}}\right)\ne0$, then
the exponential map $\exp:\mathfrak g\to G$ is not locally injective around $t_0X$.
\end{itemize}
\end{prop}
\begin{proof}
Using formula \eqref{eq:curvatura}, the Jacobi equation corresponds, via parallel transport
along $\gamma$, to the second order equation in $\mathfrak g$:
\[Y''(t)=\tfrac14\mathrm{ad}_{X}^2Y(t).\]
By Corollary \ref{cor:eigenFu}, the endomorphism $\tfrac14\mathrm{ad}_{X}^2$ has a real negative eigenvalue $\lambda$
if and only if $\tfrac12\mathrm{ad}_X$ has the purely imaginary eigenvalues $\pm i\sqrt{-\lambda}$.
Hence, part (a) of the thesis follows readily as an application of Corollary~\ref{thm:numeroconjpts}, where $A$ is the
$h$-symmetric endomorphism $\frac14\mathrm{ad}^2_{X}$ of $\mathfrak g\cong\R^n$.
Part (b) follows from Corollary~\ref{thm:cordegconjinst},  part (c) from
Corollary~\ref{thm:numeroconjpts}; the observation on the parity of the multiplicity follows from the equalities:
\begin{equation}\label{eq:dim2}
\begin{split}\mathrm{dim}_{\R}\Big(\Ker\big(&\mathrm{ad}_X^2+{4k^2\pi^2}{t_0^{-2}}\big)\Big)\\&=
\mathrm{dim}_{\mathbb C}\Big(\Ker\big(\mathrm{ad}_X-{2ik\pi}{t_0}^{-1}\big)\oplus
\Ker\big(\mathrm{ad}_X+{2ik\pi}{t_0}^{-1}\big)\Big)\\&=
2\,\mathrm{dim}_{\mathbb C}\Big(\Ker\big(\mathrm{ad}_X-{2ik\pi}{t_0}^{-1}\big)\Big).
\end{split}
\end{equation}
Here the first equality follows from Lemma \ref{lem:seriediA} and  the second one from the
fact that the two involved Kernels are conjugate spaces.
Part (d) follows from Corollary~\ref{thm:maslovindex}.
Part (e) is an application of a result on bifurcation of semi-Riemannian geodesics,
that can be found in \cite{PicPorTau}.
\end{proof}

When $\mathrm{dim}(G)\le 5$ or, more generally, when the structure coefficients
of $\mathfrak g$ satisfy the relations \eqref{eq:identities}, the statement of
Proposition~\ref{thm:gruppi} can be improved as follows:
\begin{prop}\label{thm:gruppi2}
Under the hypotheses of Proposition~\ref{thm:gruppi}, let $X_1,\ldots,X_n$
be an $h$-orthonormal basis of $\mathfrak g$ and set $\gamma(t)=\exp(tX_n)$,
for all $t\in\R$. If the structure coefficients
satisfy the relations \eqref{eq:identities}, then there are conjugate points
along $\gamma$ if and only if $\sum\limits_{i<j}\epsilon_i\epsilon_j a_{ij}^2>0$, where
the $\epsilon_{k}$'s and $a_{rs}$'s are defined as in the statement of Lemma~\ref{thm:lemmatecnico}.
In this case:
\begin{itemize}
\item[(a)] every conjugate point along $\gamma$ is nondegenerate;
\item[(b)] every conjugate point has multiplicity equal to $2$;
\item[(c)] all conjugate points along $\gamma$ give the same contribution to
the Maslov index, which is equal to $\pm2$;
\item[(d)] if $\gamma(t_0)$ is conjugate, then $\exp$ is not locally injective
around $t_0X_n$.
\end{itemize}
In particular, by (c), if $h$ has index $1$, i.e., $(G,h)$ is a Lorentzian group,
then the contribution to the Maslov index equals its multiplicity.
\end{prop}
\begin{proof}
It follows readily from Lemma~\ref{thm:lemmatecnico} and Proposition~\ref{thm:gruppi}, observing that,
for $\lambda\in i\R\setminus\{0\}$, the real generalized eigenspace $\Ker\left(\mathrm{ad}^2_{X_n}+\lambda^2\right)^n$
is the real part of the direct sum of the complex generalized eigenspaces
$\Ker\left(\mathrm{ad}_{X_n}+\lambda\right)^n$ and $\Ker\left(\mathrm{ad}_{X_n}-\lambda\right)^n$ (observe that this fact follows from Lemma \ref{lem:seriediA}).
Note that, from \eqref{eq:polcarat}, the algebraic multiplicity of each non zero eigenvalue of
$\mathrm{ad}_{X_n}$ is equal to $1$; from this observation and from part (b) of Proposition~\ref{thm:gruppi}
we obtain a proof of part (a). Moreover, it follows from Proposition~\ref{thm:gruppi} that the multiplicity
of each negative eigenvalue of $\mathrm{ad}^2_{X_n}$ is equal to $2$, which proves part (b).
Part (c) follows from the last statement in Lemma~\ref{thm:lemmatecnico} and  part (d) of Proposition~\ref{thm:gruppi}.
Finally, (d) follows
from part (e) in Proposition~\ref{thm:gruppi}.
\end{proof}

\end{subsection}
\subsection{A final remark}
The results presented in this paper are in striking contrast with
several assertions made in a preprint recently appeared \cite{Por}, where
the author attempts a calculation of the Maslov index for an
autonomous linear Hamiltonian system\footnote{%
a ``constant symplectic system'' in our terminology} using a
normal reduction for the coefficient matrix. According to what
claimed in \cite{Por}, the number, the distribution and the
contribution to the Maslov index of the conjugate instants would
depend also on the non real complex eigenvalues of the curvature
tensor (compare with Lemma~\ref{thm:existenceconjpts},
Proposition~\ref{thm:numeroconjpts} and
Corollary~\ref{thm:maslovindex}), while the lack of
diagonalizability and the role of the generalized eigenspace of
the real negative eigenvalues of the curvature tensor has not been
recognized. It is also claimed in \cite{Por} that the conjugate instants
of an arbitrary constant symplectic system do not accumulate
at the initial instant (compare with Lemma~\ref{thm:conjinstdiscr}).

Arbitrary symplectic changes of coordinates in $\R^{2n}$ needed for
the normal reduction employed in \cite{Por}
do not preserve the conjugate instants of the system
nor the base Lagrangian subspace $L_0=\{0\}\oplus\R^n$, so that a
suitable correction term has to be computed in order to get
the correct formulas. Several mistakes in such computation
and also in other parts of the preprint have led the author
to incorrect conclusions throughout.

The author of \cite{Por} has indicated the false address of the
University of S\~ao Paulo, Brazil, as his own institution; it must
be observed that A. Portaluri has no affiliation whatsoever with
such institution. The ideas, the methods and the results contained
in \cite{Por} are entire responsibility of the author's true
institution, which is the Politecnico di Torino, Italy.

\end{section}

\end{document}